\documentclass[11pt]{article}
\usepackage{enumerate}
\usepackage{lineno,hyperref}
\usepackage{amsfonts}
\usepackage{amssymb}
\usepackage{mathrsfs}
\usepackage{srcltx}
\usepackage{amsmath}
\usepackage{amsthm}
\usepackage{amstext}
\usepackage{amsopn}
\usepackage{graphicx}
\modulolinenumbers[4]

\newtheorem{theorem}{Theorem}[section]
\newtheorem{lemma}[theorem]{Lemma}
\newtheorem{proposition}[theorem]{Proposition}
\newtheorem{corollary}[theorem]{Corollary}

\theoremstyle{definition}

\theoremstyle{remark}
\newtheorem{remark}[theorem]{Remark}

\numberwithin{equation}{section}

\newcommand{\ba}{\begin{array}}
\newcommand{\ea}{\end{array}}

\newcommand{\R}{\mathbb{R}}

\newcommand{\N}{\mathbb{N}}

\allowdisplaybreaks[4]

\begin{document}
 \date{}
 \title{ Berestycki-Lions conditions on ground state solutions for a Nonlinear Schr\"odinger equation with variable potentials
\footnote{This work is partially supported by the National Natural Science Foundation of China (No: 11571370).}}
\author{ Xianhua Tang\footnote{Corresponding author.}\ and Sitong Chen\\
        {\small School of Mathematics and Statistics, Central South University,}\\
        {\small Changsha, Hunan 410083, P.R.China }\\
        {\small E-mail: tangxh@mail.csu.edu.cn(X.H. Tang)}\\
        {\small E-mail: mathsitongchen@163.com (S.T. Chen)}}
\maketitle
\begin{abstract}
 This paper is dedicated to studying the nonlinear Schr\"odinger equations of the form
 \begin{equation*}\label{KE}
 \left\{
   \begin{array}{ll}
     -\triangle u+V(x)u=f(u), & x\in \R^N; \\
     u\in H^1(\R^N),
   \end{array}
 \right.
 \end{equation*}
 where $V\in \mathcal{C}^1(\R^N, [0, \infty))$ satisfies some weak assumptions, and $f\in \mathcal{C}(\R, \R)$ satisfies the general
 Berestycki-Lions assumptions. By introducing some new tricks, we prove that the above problem admits a ground state solution
 of Poho\u zaev type and a least energy solution. These results generalize and improve some ones in [L. Jeanjean, K. Tanka,
 Indiana Univ. Math. J. 54 (2005), 443-464], [L. Jeanjean, K. Tanka, Proc. Amer. Math. Soc. 131 (2003) 2399-2408],
 [H. Berestycki, P.L. Lions,  Arch. Rational Mech. Anal. 82 (1983) 313-345] and some other related literature. In particular,
 our assumptions are ``almost" necessary when $V(x)\equiv V_{\infty}>0$, moreover, our approach could be useful
 for the study of other problems where radial symmetry of bounded sequence either fails or is not readily available,
 or where the ground state solutions of the problem at infinity are not sign definite.


\end{abstract}

{\section{Introduction}}
 \setcounter{equation}{0}

  In this paper, we consider the nonlinear Schr\"odinger equations of the form:
 \begin{equation}\label{SE}
 \left\{
   \begin{array}{ll}
     -\triangle u+V(x)u=f(u), & x\in \R^N; \\
     u\in H^1(\R^N),
   \end{array}
 \right.
 \end{equation}
 where $N\ge 3$, $V: \R^N\rightarrow \R$ and $f: \R \rightarrow \R$ satisfy the following basic assumptions:

 \begin{itemize}
 \item[(V1)] $V\in \mathcal{C}(\R^N, [0, \infty))$;

 \item[(V2)]  $V(x)\le V_{\infty}:=\lim_{|y|\to \infty}V(y)$ for all $x\in \R^N$;

 \item[(F1)]  $f\in \mathcal{C}(\R, \R)$ and there exists a constant $\mathcal{C}_0>0$ such that
 $$
   |f(t)|\le \mathcal{C}_0\left(1+|t|^{2^*-1}\right), \ \ \ \ \forall \ t\in \R;
 $$

 \item[(F2)]  $f(t)=o(t)$ as $t\to 0$ and $|f(t)|=o\left(|t|^{(N+2)/(N-2)}\right)$ as $|t|\to +\infty$.
 \end{itemize}

 \par
   Clearly, under (V1), (V2), (F1) and (F2), the weak solutions of \eqref{SE} correspond to the critical points of
 the energy functional defined in $H^1(\R^N)$ by
 \begin{equation}\label{IU}
   I(u) = \frac{1}{2}\int_{\R^N}\left[|\nabla u|^2+V(x)u^2\right]\mathrm{d}x-\int_{\R^N}F(u)\mathrm{d}x,
 \end{equation}
 where $F(t):=\int_{0}^{t}f(s)\mathrm{d}s$.

 \vskip2mm
 \par
   If the potential $V (x)\equiv V_{\infty}$, then \eqref{SE} reduces to the following autonomous form:
 \begin{equation}\label{SE1}
 \left\{
   \begin{array}{ll}
     -\triangle u+V_{\infty}u=f(u), & x\in \R^N; \\
     u\in H^1(\R^N),
   \end{array}
 \right.
 \end{equation}
 its energy functional is as follows:
 \begin{equation}\label{Ii}
   I^{\infty}(u)=\frac{1}{2}\int_{{\R}^N}\left(|\nabla u|^2+V_{\infty}u^2\right)\mathrm{d}x
            -\int_{{\R}^N}F(u)\mathrm{d}x.
 \end{equation}

 \par
    It is well known that every solution $u(x)$ of \eqref{SE1} satisfies the following Poho\u zaev type identity \cite{Po}:
 \begin{equation}\label{Pi}
  \mathcal{P}^{\infty}(u):=\frac{N-2}{2}\|\nabla u\|_2^2+\frac{NV_{\infty}}{2}\|u\|_2^2 -N\int_{{\R}^N}F(u)\mathrm{d}x=0.
 \end{equation}
 Let
 \begin{equation}\label{Mi}
   \mathcal{M}^{\infty}:=\left\{u\in H^1(\R^N)\setminus \{0\}: \mathcal{P}^{\infty}(u)=0\right\}.
 \end{equation}
 Berestycki-Lions \cite{BL} proved that \eqref{SE1} has a radially symmetric positive solution provided $f$ satisfies (F1), (F2) and the
 following two assumptions:

 \begin{itemize}
   \item[(F0)] $f$ is odd;
   \item[(F3)] there exists $s_0>0$ such that $F(s_0)>\frac{1}{2}V_{\infty}s_0^2$.
 \end{itemize}

  \par
   To prove the above result, Berestycki-Lions \cite{BL} considered the following constrained
 minimization problem
 \begin{equation}\label{min}
   \min\left\{\|\nabla u\|_2^2 : u\in\mathcal{S} \right\},
 \end{equation}
 where
 \begin{equation}\label{mS}
   \mathcal{S} = \left\{u\in H^1(\R^N) : \int_{\R^N}\left[F(u)-\frac{1}{2}V_{\infty}u^2\right]\mathrm{d}x=1\right\};
 \end{equation}
 they first showed that by the P\'olya-Szeg\"o inequality for the Schwarz symmetrization, the minimum can be taken on radial
 and radially nonincreasing functions. Then they showed the existence of a minimizer $\hat{w} \in H^1(\R^N)$ by the direct method of the
 calculus of variations. With the Lagrange multiplier Theorem, they concluded that $\bar{u}(x):=\hat{w}(x/t_{\hat{w}})$ with
 $t_{\hat{w}}=\sqrt{\frac{N-2}{2N}}\|\nabla \hat{w}\|_2$ is a least energy solution of \eqref{SE1}. By noting the one-to-one
 correspondence between $\mathcal{S}$ and $\mathcal{M}^{\infty}$, Jeanjean-Tanaka \cite{JT} proved that $\bar{u}$ is also a ground state
 solution of Poho\u zaev type for \eqref{SE1}, i.e. $\bar{u}\in \mathcal{M}^{\infty}$ and satisfies
 \begin{equation}\label{Gs}
   I^{\infty}(\bar{u})=\inf_{\mathcal{M}^{\infty}}I^{\infty}.
 \end{equation}
 By using a different way, Shatah \cite{Sh} showed that there exists $\tilde{u}\in \mathcal{M}_r^{\infty}$ such that
 \begin{equation}\label{Gsr}
   I^{\infty}(\tilde{u})=\inf_{\mathcal{M}_r^{\infty}}I^{\infty},
 \end{equation}
 where
 $$
   \mathcal{M}_r^{\infty}:=\left\{u\in H_r^1(\R^N)\setminus \{0\}: \mathcal{P}^{\infty}(u)=0\right\}
 $$
 and
 $$
   H_r^1(\R^N)=\left\{u\in H^1(\R^N) : u \mbox{ is radially symmetric function on}\ \R^N\right\}.
 $$

 \par
   Obviously, (F1)-(F3) are satisfied by a very wide class of nonlinearities. In particular only conditions on $f(t)$ near $0$,
 $\infty$ and the point $s_0$ are required. Moreover, in view of \cite[2.2]{BL}, (F1) is ``almost" necessary, and (F2) and (F3) are
 necessary for the existence of a nontrivial solution of problem \eqref{SE1}.

 \par
   When $V(x)\not\equiv V_{\infty}$, the approach used in \cite{BL} does not work any more for
 nonautonomous equation \eqref{SE}, since the Schwarz symmetrization can only be applied to autonomous problems. In a different way,
 Rabinowitz \cite{Ra} proved that \eqref{SE} has a nontrivial solution if $V$ satisfies (V1) and (V2) and $f$ does
 (F1), (F2), the Nehari monotonic condition:

 \begin{itemize}
 \item[(Ne)] $f(t)/|t|$ is strictly increasing on $(-\infty, 0)\cup (0, \infty)$;
 \end{itemize}

 \vskip2mm
 \noindent
 and the global growth Ambrosetti-Rabinowitz condition:

 \begin{itemize}
 \item[(AR)] there exists $\mu>2$ such that $f(t)t\ge \mu F(t)> 0, \ \forall \ t\in \R\setminus \{0\}$.
 \end{itemize}
 (Ne) and (AR) are used to recover the compactness and to get the boundedness of Palais-Smale sequences, respectively. By means of Jeanjean's monotonicity trick, developed in \cite{Je}, which is a generalization of the Struwe's one (see \cite{St}), consisting in a
 suitable approximating method, Jeanjean and Tanaka \cite{JT2} derived an existence result using two weaker conditions instead of (Ne) and (AR).
 More precisely, Jeanjean and Tanaka proved that \eqref{SE} has a least energy solution if $f$ satisfies (F1), (F2) and the following
 nonnegativity condition:
 \begin{itemize}
 \item[(NG)]  $f(t)\ge 0$ for $t\ge 0$;
 \end{itemize}
 and the superlinear growth condition:
 \begin{itemize}
 \item[(SL)]  $\lim_{t\to+\infty}\frac{f(t)}{t}=\infty$;
 \end{itemize}
 and $V$ does (V1), (V2) and the decay condition:
 \begin{itemize}
 \item[(Vd)] $V\in \mathcal{C}^1(\R^N, \R)$ and there exists $\varphi\in L^2(\R^N)\cap W^{1,\infty}(\R^N)$ such that
 $$
   |\nabla V(x)||x|\le [\varphi(x)]^2, \ \ \ \ \forall \ x\in \R^N.
 $$
 \end{itemize}

 \par
   Clearly, (NG) and (SL) are stronger than (F3), moreover, (Vd) puts  relatively strict constrains on the decay of $|\nabla V(x)|$.
 For example, $V(x)=a-\frac{b}{1+|x|^{\alpha}}$ does not satisfy (Vd) for $a, b>0$ and $0<\alpha\le N$.

 \par
   Motivated by \cite{BL,JT2,Sh, Ta, TC1, TC2}, we shall develop a more direct approach (the least energy squeeze approach)
 to show that \eqref{SE} has a solution $\bar{u}\in \mathcal{M}$ such that $I(\bar{u})=\inf_{\mathcal{M}}I$ under (F1)-(F3), (V1),
 (V2) and an additional decay condition on $V$:
 \begin{itemize}
 \item[(V3)] $V\in \mathcal{C}^1(\R^N, \R)$ and there exists $\theta\in [0, 1)$ such that
 \begin{eqnarray*}
   N\left[V(x)-V(tx)\right]+\left[\nabla V(x)\cdot x-\nabla V(tx)\cdot (tx)\right]
    +\frac{(N-2)^3\theta }{4t^{2}|x|^2}\left(t^{2}-1\right) \
    \left\{\begin{array}{ll}
      \ge 0, & t\ge 1,\\
      \le 0, & 0<t<1;
      \end{array}\right. \ \
 \end{eqnarray*}
 \end{itemize}
 where
 \begin{equation}\label{Ne}
   \mathcal{M}:= \left\{u\in H^1(\R^N)\setminus \{0\} : \mathcal{P}(u)=0\right\}
 \end{equation}
 and
 \begin{equation}\label{Jv}
   \mathcal{P}(u) :=  \frac{N-2}{2}\|\nabla u\|_2^2+\frac{1}{2}\int_{{\R}^N}[NV(x)+\nabla V(x)\cdot x]u^2\mathrm{d}x
    -N\int_{{\R}^N}F(u)\mathrm{d}x
 \end{equation}
 is the Poho\u zaev functional associated with \eqref{SE} (see \cite{JT2}).

 \par
   To prove the above conclusion, we shall divide our arguments into three steps: i). Choosing a minimizing sequence
 $\{u_n\}$ of $I$ on $\mathcal{M}$, which satisfies
 \begin{equation}\label{A01}
    I(u_n)\rightarrow m:=\inf_{\mathcal{M}}I, \ \ \ \ \ \mathcal{P}(u_n)= 0.
 \end{equation}
 Then showing that $\{u_n\}$ is bounded in $H^1(\R^N)$. With a concentration-compactness argument, showing that
 $\{u_n\}$ converges to some $\bar{u}\in H^1(\R^N)\setminus \{0\}$ up to translations and extraction of a subsequence.
 ii). Showing that $\bar{u}\in \mathcal{M}$ and $I(\bar{u})=\inf_{\mathcal{M}}I$. iii). Showing that $\bar{u}$ is a
 critical point of $I$. Of them, Step ii) is the most difficult due to lack of global compactness and adequate
 information on $I'(u_n)$. Since \eqref{SE} is nonautonomous, the radial compactness does not work for $\mathcal{M}$.
 To overcome this difficulty, we establish a crucial inequality related to $I(u)$, $I(u_t)$ and $\mathcal{P}(u)$
 (the IIP inequality in short, see Lemma \ref{lem 2.2}), where $u_t(x)=u(x/t)$, it plays an important role in many places
 of this paper. With the help of the IIP inequality, we then can complete Step ii) by using Lions' concentration compactness,
 the least energy squeeze approach and some subtle analysis.
 In particular, we only use Lions' concentration compactness in our arguments,
 the radial and other compactness are not required, see the proofs of Lemmas \ref{lem 2.12} and \ref{lem 3.2}. Moreover, such an approach could be
 useful for the study of other problems where radial symmetry of bounded sequence either fails or is not readily available.
 In Step iii), usually, one uses the Lagrange multipliers Theorem to show that the minimizer $\bar{u}$ is a critical point of $I$,
 but it is impossible to verify $\mathcal{P}'(u)\ne 0$ for all $u\in\mathcal{M}$ under (V1)-(V3) and (F1)-(F3). To overcome
 this difficulty, we employ the combination of the IIP inequality, a deformation lemma and the degree theory, see Lemma \ref{lem 2.13}.

 \begin{remark}\label{rem1.1} There are indeed functions which satisfy (V1)-(V3). An example is given by $V(x)
 =V_1-\frac{A}{|x|^2+1}$, where $V_1\ge A$ and $0<A<(N-2)^3/2(2N+1)$ are two positive constants.
 \end{remark}

 \vskip4mm
 \par
   We are now in a position to state the main results on ground state solutions of Poho\u zaev type.

 \begin{theorem}\label{thm1.2} Assume that $V$ and $f$ satisfy {\rm(V1)-(V3)} and {\rm(F1)-(F3)}. Then
 problem \eqref{SE} has a solution $\bar{u}\in H^1(\R^N)$ such that $I(\bar{u})=\inf_{\mathcal{M}}I
 =\inf_{u\in \Lambda}\max_{t > 0}I(u_t)>0$, where
 $$
   u_t(x):=u(x/t) \ \ \mbox{and} \ \
  \Lambda：=\left\{u\in H^1(\R^N) : \int_{\R^N}\left[\frac{1}{2}V_{\infty}u^2-F(u)\right]\mathrm{d}x<0\right\}.
 $$
 \end{theorem}

 \vskip4mm
 \par
  As a consequence of Theorem \ref{thm1.2}, we can prove the following theorem.

  \begin{theorem}\label{thm1.3} Assume that $f$ satisfies {\rm(F1)-(F3)}. Then problem \eqref{SE1} has a solution
 $\bar{u}\in H^1(\R^N)$ such that $I^{\infty}(\bar{u})=\inf_{\mathcal{M}^{\infty}}I^{\infty}
 =\inf_{u\in \Lambda}\max_{t > 0}I^{\infty}(u_t)>0$.
 \end{theorem}

 \begin{remark}\label{rem1.4} We point out that, as a consequence of Theorem \ref{thm1.2}, the least energy value $m:=\inf_{\mathcal{M}}I$
 has a minimax  characterization $ m=\inf_{u\in \Lambda}\max_{t > 0}I(u_t)$ which is much simpler than the usual characterizations
 related to the Mountain Pass level.
 \end{remark}

 \vskip4mm
 \par
   In the second part of the paper, we are interested in the existence of the least energy solutions for \eqref{SE} under
 (F1)-(F3). In this case, we can replace (V3) by the following weaker decay assumption on $\nabla V$:
 \begin{itemize}
 \item[(V4)] $V\in \mathcal{C}^1(\R^N, \R)$ and there exists $\theta\in [0, 1)$ such that
 $$
   \nabla V(x)\cdot x\le  \frac{(N-2)^2\theta}{2|x|^2}, \ \ \ \ \forall \ x\in \R^N\setminus \{0\}.
 $$
 \end{itemize}

 \par
   As in Jeanjean and Tanaka \cite{JT}, for $\lambda\in [1/2, 1]$ we consider a family of functionals $I_{\lambda}
  : H^1(\R^N) \rightarrow \R$ defined by
 \begin{equation}\label{Ilu}
   I_{\lambda}(u)=\frac{1}{2}\int_{\R^N}\left(|\nabla u|^2+V(x)u^2\right)\mathrm{d}x
            -\lambda\int_{\R^N}F(u)\mathrm{d}x.
 \end{equation}
 These functionals have a Mountain Pass geometry, and denoting the corresponding Mountain Pass levels by $c_{\lambda}$.
 Corresponding to \eqref{Ilu}, we also let
 \begin{equation}\label{Il}
   I_{\lambda}^{\infty}(u)=\frac{1}{2}\int_{\R^N}\left(|\nabla u|^2+V_{\infty}u^2\right)\mathrm{d}x
       -\lambda\int_{\R^N}F(u)\mathrm{d}x.
 \end{equation}
 By Theorem \ref{thm1.3}, for every $\lambda\in [1/2, 1]$, there exists a minimizer $u_{\lambda}^{\infty}$ of $I_{\lambda}^{\infty}$
 on $\mathcal{M}_{\lambda}^{\infty}$, where
 \begin{eqnarray}\label{Ml}
   \mathcal{M}_{\lambda}^{\infty}:=\left\{u\in H^1(\R^N)\setminus \{0\}: \mathcal{P}_{\lambda}^{\infty}(u)=0\right\}
 \end{eqnarray}
 and
 \begin{eqnarray}\label{JlL}
  \mathcal{P}_{\lambda}^{\infty}(u)
    &  =  & \frac{N-2}{2}\|\nabla u\|_2^2+NV_{\infty}\|u\|_2^2-N\lambda\int_{\R^N}F(u)\mathrm{d}x.
 \end{eqnarray}
 Let
 $$
   A(u)=\frac{1}{2}\int_{\R^N}\left(|\nabla u|^2+V(x)u^2\right)\mathrm{d}x, \ \ \ \
   B(u)=\frac{1}{2}\int_{\R^N}F(u)\mathrm{d}x.
 $$
 Then $I_{\lambda}(u)=A(u)-\lambda B(u)$. Since $B(u)$ is not sign definite, it prevents us from employing Jeanjean's
 monotonicity trick \cite{Je}. Thanks to the work of Jeanjean-Toland \cite{JTo}, $I_{\lambda}$ still has a bounded (PS)-sequence
 $\{u_n(\lambda)\} \subset H^1(\R^N)$ at level $c_{\lambda}$ for almost every $\lambda\in [1/2,1]$. Different from the
 arguments in the existing literature, by means of $u_1^{\infty}$ and the IIP inequality, we can find $\bar{\lambda}\in [1/2, 1)$
 and then directly prove the following crucial inequality
 \begin{eqnarray}\label{cm1}
   c_{\lambda}<m_{\lambda}^{\infty}:=\inf_{\mathcal{M}_{\lambda}^{\infty}}I_{\lambda}^{\infty},\ \ \ \ \lambda\in (\bar{\lambda}, 1],
 \end{eqnarray}
 which is used to recover the compactness to (PS)-sequence $\{u_n(\lambda)\}$, see Lemmas \ref{lem 4.5} and \ref{lem 4.7}.
 In particular, it is not required any information on sign of $u_1^{\infty}$ in our arguments. Applying \eqref{cm1} and
 a precise decomposition of bounded (PS)-sequence in \cite{JT}, we can get a nontrivial critical point $u_{\lambda}$ of
 $I_{\lambda}$ which possesses energy $c_{\lambda}$  for almost every $\lambda\in (\bar{\lambda}, 1]$. Finally,
 with a Poho\u zaev identity  we proved that \eqref{SE} admits a least energy solution under (V1), (V2), (V4) and (F1)-(F3).
 More precisely, we have the following theorem.

 \par
 \begin{theorem}\label{thm1.5} Assume that $V$ and $f$ satisfy {\rm(V1), (V2), (V4)}  and {\rm(F1)-(F3)}. Then problem \eqref{SE}
 has a least energy solution.
 \end{theorem}

 \begin{remark}\label{rem1.6} Relative to (Vd), there seem to be more functions satisfying (V4). For example, it is
 easy to verify that $V(x)=a-\frac{b}{1+|x|^{\alpha}}$ satisfies (V4) for $\alpha\ge 2$, $a>0$ and
 $$
    4b<(N-2)^2 \ \ \mbox{if} \ \alpha=2; \ \ \ \ \frac{(\alpha-2)^{(\alpha-2)/\alpha}(\alpha+2)^{(\alpha+2)/\alpha}}{2\alpha}b<(N-2)^2,
     \ \ \mbox{if} \ \alpha>2.
 $$
 However, it does not satisfy (Vd)  when $2\le \alpha \le N$.
 \end{remark}

 \vskip4mm
 \par
   Applying Theorem \ref{thm1.5} to the following perturbed problem:
 \begin{equation}\label{KE8}
 \left\{
   \begin{array}{ll}
     -\triangle u+[V_{\infty}-\varepsilon h(x)]u=f(u), & x\in \R^N; \\
     u\in H^1(\R^N),
   \end{array}
 \right.
 \end{equation}
 where $V_{\infty}$ is a positive constant and the function $h \in \mathcal{C}^1(\R^N, \R)$ verifies:

 \vskip2mm
 \noindent
 (H1)\ \ $h(x)\ge 0$ for all $x\in \R^N$ and $\lim_{|x|\to\infty}h(x)=0$;

 \vskip2mm
 \noindent
 (H2)\ \ $\sup_{x\in \R^N}\left[-|x|^2(\nabla h(x), x)\right]<\infty$.

 \vskip2mm
 \noindent
 Then we have the following corollary.

\begin{corollary}\label{cor1.7}
 Assume that $h$ and $f$ satisfy {\rm(H1), (H2)} and {\rm(F1)-(F3)}. Then there
 exists a constant $\varepsilon_0>0$ such that problem \eqref{KE8} has a least energy solution $\bar{u}_{\varepsilon}\in
 H^1(\R^N)\setminus \{0\}$ for all $0<\varepsilon\le \varepsilon_0$.
\end{corollary}

 \vskip4mm
 \par
    Classically, in order to show the existence of solutions for \eqref{SE}, one compares the critical level of $I$ with the one of
 $I^{\infty}$ (i.e. the energy functional correspondings to the problem at infinity). To this end, it is necessary to establish
 a strict inequality similar to
 $$
   \max_{t\in [0,1]}I(\gamma_0(t))<\inf\left\{I^{\infty}(u): u\in H^1(\R^N)\setminus \{0\}
   \ \mbox{is a solution of}\ \eqref{SE1}\right\}
 $$
 for some path $\gamma_0\in \mathcal{C}([0,1], H^1(\R^N))$. Clearly, $\gamma_0(t)>0$ is a natural requirement under (V1) and (V2).
 But we only need $\gamma_0(t)\ne 0$ in our arguments. Therefore, our approach could be useful for the study of other problems where paths or the ground state solutions of the problem at infinity are not sign definite.

 \vskip4mm
 \par
   Throughout the paper we make use of the following notations:

 \vskip4mm
 \par
     $\spadesuit$ \ $H^1(\R^N)$ denotes the usual Sobolev space equipped with the inner product and norm
 $$
   (u, v)=\int_{\R^N}(\nabla u\cdot \nabla v+uv)\mathrm{d}x, \ \ \|u\|=(u, u)^{1/2},
     \ \ \forall \ u,v\in H^1(\R^N);
 $$

 \par
     $\spadesuit$ \ $L^s(\R^N) (1\le s< \infty)$  denotes the Lebesgue space with the norm $\|u\|_s
 =\left(\int_{\R^N}|u|^s\mathrm{d}x\right)^{1/s}$;

 \par
     $\spadesuit$ \ For any $u\in H^1(\R^N)\setminus \{0\}$, $u_t(x):=u(x/t)$ for $t>0$;

 \par
     $\spadesuit$ \ For any $x\in \R^N$ and $r>0$, $B_r(x):=\{y\in \R^N: |y-x|<r \}$;

 \par
     $\spadesuit$ \ $C_1, C_2,\cdots$ denote positive constants possibly different in different places.

 \vskip4mm
   The rest of the paper is organized as follows. In Section 2, we give some preliminaries, and give the proof
 of Theorem \ref{thm1.3}. In Section 3, we complete the proof of Theorem \ref{thm1.2}. Section 4 is devoted to finding a least energy
 solution for \eqref{SE} and Theorem \ref{thm1.5} will be proved in this section.

 \vskip6mm

 {\section{Ground state solutions for the ``limited problem"}}
 \setcounter{equation}{0}

 \vskip2mm
 \par
   In this section, we give the proof of Theorem \ref{thm1.3}. To this end, we first give some useful lemmas. Since $V(x)\equiv V_{\infty}$
 satisfies (V1)-(V3), thus all conclusions on $I$ are also true for $I^{\infty}$. By a simple calculation, we can verify Lemma \ref{lem 2.1}.

 \vskip2mm
 \par
 \begin{lemma}\label{lem 2.1}
 The following inequality holds:
 \begin{equation}\label{B20}
  \mathfrak{g}(t):=2-Nt^{N-2}+(N-2)t^N > \mathfrak{g}(1)=0,   \ \ \ \ \forall \ t\in [0, 1)\cup(1, +\infty).
 \end{equation}
  Moreover {\rm(V3)} implies the following inequality holds:
 \begin{eqnarray}\label{V3}
    &     & Nt^N\left[V(x)-V(tx)\right]
             +\left(t^{N}-1\right)\nabla V(x)\cdot x \nonumber\\
    & \ge & -\frac{(N-2)^2\theta\left[2-Nt^{N-2}+(N-2)t^{N}\right] }{4|x|^2},
              \ \ \ \ \forall \ t\ge 0, \ \ x\in \R^N\setminus \{0\}.
 \end{eqnarray}
\end{lemma}

 \begin{lemma}\label{lem 2.2}
 Assume that {\rm(V1), (V3), (F1)} and {\rm(F2)} hold. Then
 \begin{eqnarray}\label{G31}
   I(u) & \ge & I(u_t)+\frac{1-t^N}{N}\mathcal{P}(u)+\frac{(1-\theta)\left[2-Nt^{N-2}+(N-2)t^N\right]}{2N}\|\nabla u\|_2^2,\nonumber\\
        &     & \ \ \ \ \forall \ u\in H^1(\R^N), \ \ t > 0.
 \end{eqnarray}
\end{lemma}

\begin{proof}  According to Hardy inequality, we have
 \begin{equation}\label{G38}
  \|\nabla u\|_2^2 \ge \frac{(N-2)^2}{4}\int_{\R^N}\frac{u^2}{|x|^2}\mathrm{d}x,  \ \ \forall \ u\in H^1(\R^N).
 \end{equation}
 Note that
 \begin{equation}\label{G33}
   I(u_t) = \frac{t^{N-2}}{2}\|\nabla u\|_2^2+\frac{t^N}{2}\int_{{\R}^N}V(tx)u^2\mathrm{d}x
             -t^N\int_{{\R}^N}F(u)\mathrm{d}x.
 \end{equation}
 Thus, by \eqref{IU}, \eqref{Jv}, \eqref{B20}, \eqref{V3}, \eqref{G38} and \eqref{G33}, one has
 \begin{eqnarray*}
   &     & I(u)-I(u_t)\\
   &  =  & \frac{1-t^{N-2}}{2}\|\nabla u\|_2^2+\frac{1}{2}\int_{{\R}^N}\left[V(x)-t^{N}V(tx)\right]u^2\mathrm{d}x
             -\left(1-t^N\right)\int_{\R^N}F(u)\mathrm{d}x\\
   &  =  & \frac{1-t^N}{N}\left\{\frac{N-2}{2}\|\nabla u\|_2^2+\frac{1}{2}\int_{\R^N}[NV(x)+\nabla V(x)\cdot x]u^2\mathrm{d}x
             -N\int_{{\R}^N}F(u)\mathrm{d}x\right\}\\
   &     & +\frac{2-Nt^{N-2}+(N-2)t^N}{2N}\|\nabla u\|_2^2\\
   &     &   +\frac{1}{2}\int_{\R^N}\left\{t^N[V(x)-V(tx)]
             -\frac{1-t^N}{N}\nabla V(x)\cdot x\right\}u^2\mathrm{d}x\\
   & \ge & \frac{1-t^N}{N}\mathcal{P}(u)+\frac{(1-\theta)\left[2-Nt^{N-2}+(N-2)t^N\right]}{2N}\|\nabla u\|_2^2,
             \ \ \ \ \forall \ u\in H^1(\R^N), \ \ t > 0.
 \end{eqnarray*}
 This shows that \eqref{G31} holds.
 \end{proof}

 \vskip4mm
 \par
   From Lemma \ref{lem 2.2}, we have the following two corollaries.

 \begin{corollary}\label{cor 2.3}
 Assume that {\rm(F1)} and {\rm(F2)} hold. Then
 \begin{eqnarray}\label{G35}
   I^{\infty}(u)
     & \ge & I^{\infty}(u_t)+\frac{1-t^N}{N}\mathcal{P}^{\infty}(u)+\frac{2-Nt^{N-2}+(N-2)t^N}{2N}\|\nabla u\|_2^2,\nonumber\\
     &     &    \ \ \ \ \forall \ u\in H^1(\R^N), \ \ t > 0.
 \end{eqnarray}
\end{corollary}

 \begin{corollary}\label{cor 2.4}
 Assume that {\rm(V1), (V3), (F1)} and {\rm(F2)} hold. Then for $u\in \mathcal{M}$
 \begin{equation}\label{Imax}
   I(u) = \max_{t> 0}I(u_t).
 \end{equation}
\end{corollary}

 \begin{lemma}\label{lem 2.5}
 Assume that {\rm(V1)-(V3)} hold. Then there exist two constants $\gamma_1, \gamma_2>0$
 such that
 \begin{equation}\label{G36}
   \gamma_1\|u\|^2\le (N-2)\|\nabla u\|_2^2+\int_{\R^N}\left[NV(x)+\nabla V(x)\cdot x\right]u^2\mathrm{d}x\le \gamma_2\|u\|^2,
     \ \ \forall \ u\in H^1(\R^N).
 \end{equation}
\end{lemma}

\begin{proof}  Let $t=0$ and $t \to \infty$ in \eqref{V3}, respectively, and using (V2), one has
 \begin{equation}\label{G37}
   -\frac{(N-2)^3\theta}{4|x|^2}+NV_{\infty} \le NV(x)+\nabla V(x)\cdot x\le NV_{\infty}+\frac{(N-2)^2\theta}{2|x|^2},
       \ \ \ \ \forall \ x\in \R^N\setminus \{0\}.
 \end{equation}
 Thus it follows from \eqref{G38} and \eqref{G37} that
 \begin{eqnarray}\label{G40}
   &     & (N-2)\|\nabla u\|_2^2+\int_{\R^N}\left[NV(x)+\nabla V(x)\cdot x\right]u^2\mathrm{d}x\nonumber\\
   & \le & (N-2+2\theta)\|\nabla u\|_2^2+NV_{\infty}\|u\|_2^2\nonumber\\
   & \le & [N-2+2\theta+NV_{\infty}]\|u\|^2:=\gamma_2\|u\|^2, \ \ \forall \ u\in H^1(\R^N)
 \end{eqnarray}
 and
 \begin{eqnarray}\label{G41}
   &     & (N-2)\|\nabla u\|_2^2+\int_{\R^N}\left[NV(x)+\nabla V(x)\cdot x\right]u^2\mathrm{d}x\nonumber\\
   & \ge & (1-\theta)(N-2)\|\nabla u\|_2^2+NV_{\infty}\|u\|_2^2\nonumber\\
   & \ge & \min\left\{(1-\theta)(N-2), NV_{\infty}\right\}\|u\|^2:=\gamma_1\|u\|^2, \ \ \forall \ u\in H^1(\R^N).
 \end{eqnarray}
 Both \eqref{G40} and \eqref{G41} imply that \eqref{G36} holds.
 \end{proof}

 \vskip4mm
 \par
   To show $\mathcal{M}\ne \emptyset$, we define a set $\Lambda$ as follows:
 \begin{equation}\label{La}
   \Lambda：=\left\{u\in H^1(\R^N) : \int_{\R^N}\left[\frac{1}{2}V_{\infty}u^2-F(u)\right]\mathrm{d}x<0\right\}.
 \end{equation}

 \begin{lemma}\label{lem 2.6}
 Assume that {\rm(V1)-(V3)} and {\rm(F1)-(F3)} hold. Then $\Lambda\ne\emptyset$ and
 \begin{equation}\label{La1}
   \left\{u\in H^1(\R^N)\setminus \{0\} : \mathcal{P}^{\infty}(u)\le 0 \ \mbox{or} \ \mathcal{P}(u)\le 0\right\}\subset \Lambda.
 \end{equation}
\end{lemma}

\begin{proof} In view of the proof of \cite[Theorem 2]{BL}, (F3) implies $\Lambda\ne\emptyset$. Next, we have two cases
 to distinguish:

 \vskip2mm
 \par
  1). $u\in H^1(\R^N)\setminus \{0\}$ and $\mathcal{P}^{\infty}(u)\le 0$, then \eqref{Pi} implies $u\in \Lambda$.

 \par
  2). $u\in H^1(\R^N)\setminus \{0\}$ and $\mathcal{P}(u)\le 0$, then it follows from \eqref{Jv}, \eqref{G38} and \eqref{G37} that
 \begin{eqnarray*}
   &     &  N\int_{\R^N}\left[\frac{1}{2}V_{\infty}u^2-F(u)\right]\mathrm{d}x\\
   &  =  & \mathcal{P}(u)-\frac{N-2}{2}\|\nabla u\|_2^2-\frac{N}{2}\int_{{\R}^N}\left[(V(x)-V_{\infty})
             +\frac{\nabla V(x)\cdot x}{N}\right]u^2\mathrm{d}x\\
   & \le & -\frac{N-2}{2}\|\nabla u\|_2^2+\frac{(N-2)^3\theta}{8}\int_{{\R}^N}\frac{u^2}{|x|^2}\mathrm{d}x\\
   & \le & -\frac{(1-\theta)(N-2)}{2}\|\nabla u\|_2^2 < 0,
 \end{eqnarray*}
 which implies $u\in \Lambda$.
 \end{proof}

 \begin{lemma}\label{lem 2.7}
 Assume that {\rm(V1)-(V3)} and {\rm(F1)-(F3)} hold. Then for any
 $u\in \Lambda$, there exists a unique $t_u>0$ such that $u_{t_u}\in \mathcal{M}$.
\end{lemma}

\begin{proof}  Let $u\in \Lambda$ be fixed and define a function $\zeta(t):=I(u_t)$
 on $(0, \infty)$. Clearly, by \eqref{Jv} and \eqref{G33}, we have
 \begin{eqnarray}\label{B40}
    &     & \zeta'(t)=0 \nonumber\\
    &     & \Leftrightarrow \ \ \frac{N-2}{2}t^{N-2}\|\nabla u\|_2^2+\frac{t^N}{2}\int_{{\R}^N}[NV(tx)
              +\nabla V(tx)\cdot (tx)]u^2\mathrm{d}x -Nt^N\int_{{\R}^N}F(u)\mathrm{d}x=0\nonumber\\
    &     & \ \Leftrightarrow  \ \ \mathcal{P}(u_t)=0 \ \ \Leftrightarrow  \ \ u_t\in \mathcal{M}.
 \end{eqnarray}
 It is easy to verify, using (V1), (V2), (F1), \eqref{G33} and the definition of $\Lambda$, that $\lim_{t\to 0}\zeta(t)=0$, $\zeta(t)>0$
 for $t>0$ small and $\zeta(t)<0$ for $t$ large. Therefore $\max_{t\in [0, \infty)}\zeta(t)$ is achieved at some $t_u>0$ so that
 $\zeta'(t_u)=0$ and $u_{t_u}\in \mathcal{M}$.

 \par
    Next we claim that $t_u$ is unique for any $u\in \Lambda$. In fact, for any given $u\in \Lambda$,
 let $t_1, t_2>0$ such that $u_{t_1}, u_{t_2} \in \mathcal{M}$. Then $\mathcal{P}\left(u_{t_1}\right)=\mathcal{P}\left(u_{t_2}\right)=0$.
 Jointly with \eqref{G31}, we have
 \begin{eqnarray}\label{B41}
   I\left(u_{t_1}\right)
   & \ge & I\left(u_{t_2}\right)+\frac{t_1^N-t_2^N}{Nt_1^N}\mathcal{P}\left(u_{t_1}\right)
             +\frac{(1-\theta)\left[2t_1^N-Nt_1^2t_2^{N-2}+(N-2)t_2^N\right]}{2Nt_1^N}\|\nabla u_{t_1}\|_2^2\nonumber\\
   &  =  & I\left(u_{t_2}\right)+\frac{(1-\theta)\left[2t_1^N-Nt_1^2t_2^{N-2}+(N-2)t_2^N\right]}{2Nt_1^2}\|\nabla u\|_2^2
 \end{eqnarray}
 and
 \begin{eqnarray}\label{B42}
   I\left(u_{t_2}\right)
   & \ge & I\left(u_{t_1}\right)+\frac{t_2^N-t_1^N}{Nt_2^N}\mathcal{P}\left(u_{t_2}\right)
             +\frac{(1-\theta)\left[2t_2^N-Nt_2^2t_1^{N-2}+(N-2)t_1^N\right]}{2Nt_2^N}\|\nabla u_{t_2}\|_2^2\nonumber\\
   &  =  & I\left(u_{t_1}\right)+\frac{(1-\theta)\left[2t_2^N-Nt_2^2t_1^{N-2}+(N-2)t_1^N\right]}{2Nt_2^2}\|\nabla u\|_2^2.
 \end{eqnarray}
 \eqref{B41} and \eqref{B42} imply $t_1=t_2$. Therefore, $t_u> 0$ is unique for any $u\in \Lambda$.
\end{proof}

\begin{corollary}\label{cor 2.8}
 Assume that {\rm(F1)-(F3)} hold. Then for any
 $u\in \Lambda$, there exists a unique $t_u>0$ such that $u_{t_u}\in \mathcal{M}^{\infty}$.
\end{corollary}

 \vskip4mm
 \par
    From Corollary \ref{cor 2.4}, Lemma \ref{lem 2.6} and Lemma \ref{lem 2.7}, we have $\mathcal{M}\ne \emptyset$ and the following lemma.

 \begin{lemma}\label{lem 2.9}
 Assume that {\rm(V1)-(V3)} and {\rm(F1)-(F3)} hold. Then
 $$
   \inf_{u\in \mathcal{M}}I(u)
   :=m=\inf_{u\in \Lambda}\max_{t > 0}I(u_t).
 $$
\end{lemma}

 \vskip2mm
 \par
  The following lemma is a known result which can be proved by a standard argument.

 \begin{lemma}\label{lem 2.10}
 Assume that {\rm(V1), (V2), (F1)} and {\rm(F2)} hold. If $u_n\rightharpoonup \bar{u}$ in $H^1(\R^N)$, then
 \begin{equation}\label{F60}
   I(u_n)=I(\bar{u})+I(u_n-\bar{u})+o(1)
 \end{equation}
 and
 \begin{equation}\label{F63}
   \mathcal{P}(u_n)=\mathcal{P}(\bar{u})+\mathcal{P}(u_n-\bar{u})+o(1).
 \end{equation}
\end{lemma}

\begin{lemma}\label{lem 2.11}
 Assume that {\rm(V1)-(V3)} and {\rm(F1)-(F3)} hold. Then
 \begin{enumerate}[{\rm(i)}]
  \item there exists $\rho_0>0$ such that $\|u\|\ge \rho_0, \ \forall \ u\in \mathcal{M}$;
  \item $m=\inf_{u\in \mathcal{M}}I(u)>0$.
 \end{enumerate}
\end{lemma}

\begin{proof} i). Since $\mathcal{P}(u)=0, \ \forall u\in \mathcal{M}$, by (F1), (F2), \eqref{Jv}, \eqref{G36} and Sobolev
 embedding theorem, one has
 \begin{eqnarray}\label{G62}
   \frac{\gamma_1}{2}\|u\|^2
   & \le & \frac{N-2}{2}\|\nabla u\|_2^2+\frac{1}{2}\int_{{\R}^N}[NV(x)+\nabla V(x)\cdot x]u^2\mathrm{d}x \nonumber\\
   &  =  & N\int_{{\R}^N}F(u)\mathrm{d}x \nonumber\\
   & \le & \frac{\gamma_1}{4}\|u\|^2+C_1\|u\|^{2^*},
 \end{eqnarray}
 which implies
 \begin{equation}\label{G63}
   \|u\|\ge \rho_0:=\left(\frac{\gamma_1}{4C_1}\right)^{(N-2)/4}, \ \ \ \ \forall \ u\in \mathcal{M}.
 \end{equation}

 \par
   ii). For $u\in H^1(\R^N)$, by the Sobolev inequality, one has $S\|u\|_{2^*}^2\le \|\nabla u\|_2^2$.
 By (V2), there exists $R>0$ such that $V(x)\ge \frac{V_{\infty}}{2}$ for $|x|\ge R$.  It follows from (F1) and (F2) that
 there exists $C_2>0$ such that
 \begin{equation}\label{G64}
   |F(t)|\le  \frac{1}{4}\min\left\{\frac{S}{R^2\omega_N^{2/N}},
              V_{\infty}\right\}|t|^2+C_2|t|^{2^*}, \ \ \ \ \forall \ t\in \R,
 \end{equation}
 where $\omega_N$ denote the volume of the unit ball of $\R^N$. For $u\in \mathcal{M}$, let
 $$
   t_u =\left[\frac{(N-2)S^{N/(N-2)}}{4NC_2}\right]^{1/2}\|\nabla u\|_2^{-2/(N-2)}.
 $$
 Making use of the H\"older inequality and the Sobolev embedding theorem, we get
 \begin{equation}\label{G65}
   \int_{|t_ux|< R}u^2\mathrm{d}x\le \left(\frac{\omega_N R^N}{t_u^N}\right)^{(2^*-2)/2^*}\left(\int_{|t_ux|< R}u^{2^*}\mathrm{d}x\right)^{2/2^*}
     \le \left(\frac{\omega_N R^N}{t_u^N}\right)^{2/N}S^{-1}\|\nabla u\|_2^2.
 \end{equation}
 Then from \eqref{G31}, \eqref{G33}, \eqref{G64}, \eqref{G65} and the Sobolev embedding theorem, we have
 \begin{eqnarray*}
   I(u)
    & \ge & I\left(u_{t_u}\right)\nonumber\\
    &  =  &  \frac{t_u^{N-2}}{2}\|\nabla u\|_2^2+\frac{t_u^N}{2}\int_{{\R}^N}V(t_ux)u^2\mathrm{d}x
             -t_u^N\int_{{\R}^N}F(u)\mathrm{d}x\nonumber\\
    & \ge & \frac{t_u^{N-2}}{4}\|\nabla u\|_2^2+\frac{S}{4R^2\omega_N^{2/N}}t_u^N\int_{|t_ux|< R}u^2\mathrm{d}x
              +\frac{V_{\infty}t_u^N}{4}\int_{|t_ux|\ge R}u^2\mathrm{d}x\nonumber\\
    &     &   -t_u^N\int_{\R^N}F(u)\mathrm{d}x\nonumber\\
    & \ge & \frac{t_u^{N-2}}{4}\|\nabla u\|_2^2+\frac{1}{4}\min\left\{\frac{S}{R^2\omega_N^{2/N}},
              V_{\infty}\right\}t_u^N\|u\|_2^2-t_u^N\int_{\R^N}F(u)\mathrm{d}x\nonumber\\
    & \ge & \frac{t_u^{N-2}}{4}\|\nabla u\|_2^2-C_2t_u^{N}\|u\|_{2^*}^{2^*}\nonumber\\
    & \ge & \frac{t_u^{N-2}}{4}\|\nabla u\|_2^2-C_2S^{-N/(N-2)}t_u^{N}\|\nabla u\|_2^{2N/(N-2)}\nonumber\\
    &  =  & \frac{S^{N/2}}{2^{N-1}(N-2)C_2^{(N-2)/2}}\left(\frac{N-2}{N}\right)^{N/2}, \ \ \ \  \forall \ u\in \mathcal{M}.
 \end{eqnarray*}
 This shows that $m=\inf_{u\in \mathcal{M}}I(u)>0$.
 \end{proof}

 \begin{lemma}\label{lem 2.12}
 Assume that {\rm(F1)-(F3)} hold. Then $m^{\infty}:=\inf_{u\in \mathcal{M}^{\infty}}I^{\infty}(u)$ is achieved.
\end{lemma}

\begin{proof} In view of Lemma \ref{lem 2.6}, Corollary \ref{cor 2.8} and Lemma \ref{lem 2.11}, we have $\mathcal{M}^{\infty}\ne \emptyset$
 and $m^{\infty}>0$. Let $\{u_n\}\subset \mathcal{M}^{\infty}$ be such that $I^{\infty}(u_n)\rightarrow m^{\infty}$.
 Since $\mathcal{P}^{\infty}(u_n)=0$, it follows from \eqref{G35} with $t\to 0$ that
 \begin{equation}\label{B74}
   m^{\infty}+o(1)= I^{\infty}(u_n)\ge \frac{1-\theta}{N}\|\nabla u_n\|_2^2.
 \end{equation}
 This shows that $\{\|\nabla u_n\|_2\}$ is bounded. Next, we prove that $\{\|u_n\|\}$ is also bounded.
 By (F1), (F2), \eqref{Pi} and Sobolev embedding theorem, one has
 \begin{eqnarray}\label{B75}
   \min\{N-2, NV_{\infty}\}\|u_n\|^2
   & \le & (N-2)\|\nabla u\|_2^2+NV_{\infty}\|u\|_2^2\nonumber\\
   &  =  & 2N\int_{{\R}^N}F(u)\mathrm{d}x \nonumber\\
   & \le & \frac{1}{2}\min\{N-2, NV_{\infty}\}\|u_n\|^2+C_3\|u_n\|_{2^*}^{2^*}\nonumber\\
   & \le & \frac{1}{2}\min\{N-2, NV_{\infty}\}\|u_n\|^2+C_3S^{-2^*/2}\|\nabla u_n\|_2^{2^*}.
 \end{eqnarray}
 Hence, $\{u_n\}$ is bounded in $H^1(\R^N)$. By Lions' concentration compactness principle
 \cite[Lemma 1.21]{WM}, one can easily prove that there exist $\delta>0$ and $\{y_n\}\subset \R^N$ such that
 $\int_{B_1(y_n)}|u_n|^2\mathrm{d}x> \delta/2$. Let $\hat{u}_n(x)=u_n(x+y_n)$. Then $\|\hat{u}_n\|=\|u_n\|$,
 \begin{equation}\label{B79}
   \int_{B_1(0)}|\hat{u}_n|^2\mathrm{d}x> \frac{\delta}{2}
 \end{equation}
 and
 \begin{equation}\label{D71}
    I^{\infty}(\hat{u}_n)\rightarrow m^{\infty}, \ \ \ \ \ \mathcal{P}^{\infty}(\hat{u}_n)= 0.
 \end{equation}
 Therefore, there exists $\hat{u}\in H^1(\R^N)\setminus \{0\}$ such that
 \begin{equation}\label{D72}
 \left\{
   \begin{array}{ll}
     \hat{u}_n\rightharpoonup \hat{u}, & \mbox{in} \ H^1(\R^N); \\
     \hat{u}_n\rightarrow \hat{u}, & \mbox{in} \ L_{\mathrm{loc}}^s(\R^N), \ \forall \ s\in [1, 2^*);\\
     \hat{u}_n\rightarrow \hat{u}, & \mbox{a.e. on} \ \R^N.
   \end{array}
 \right.
 \end{equation}
 Let $w_n=\hat{u}_n-\hat{u}$. Then \eqref{D72} and Lemma \ref{lem 2.10} yield
 \begin{equation}\label{D73}
    I^{\infty}(\hat{u}_n) = I^{\infty}(\hat{u})+I^{\infty}(w_n)+o(1).
 \end{equation}
 and
 \begin{equation}\label{D74}
    \mathcal{P}^{\infty}(\hat{u}_n) = \mathcal{P}^{\infty}(\hat{u})+\mathcal{P}^{\infty}(w_n)+o(1).
 \end{equation}
 From \eqref{Ii}, \eqref{Pi}, \eqref{D71}, \eqref{D73} and \eqref{D74}, one has
 \begin{equation}\label{D75}
    \frac{1}{N}\|\nabla w_n\|_2^2=m^{\infty}-\frac{1}{N}\|\nabla \hat{u}\|_2^2+o(1),
      \ \ \ \ \mathcal{P}^{\infty}(w_n) = -\mathcal{P}^{\infty}(\hat{u})+o(1).
 \end{equation}
 If there exists a subsequence $\{w_{n_i}\}$ of $\{w_n\}$ such that $w_{n_i}=0$, then going to this subsequence, we have
 \begin{equation}\label{D76}
    I^{\infty}(\hat{u})=m^{\infty}, \ \ \ \ \mathcal{P}^{\infty}(\hat{u})=0,
 \end{equation}
 which implies the conclusion of Lemma \ref{lem 2.12} holds. Next, we assume that $w_n\ne 0$. We claim that $\mathcal{P}^{\infty}(\hat{u})\le 0$.
 Otherwise, if $\mathcal{P}^{\infty}(\hat{u})>0$, then \eqref{D75} implies $\mathcal{P}^{\infty}(w_n) < 0$ for large $n$.
 In view of Lemma \ref{lem 2.6} and Corollary \ref{cor 2.8}, there exists $t_n>0$ such that $(w_n)_{t_n}\in \mathcal{M}^{\infty}$. From \eqref{Ii}, \eqref{Pi}, \eqref{G35}
 and \eqref{D75}, we obtain
 \begin{eqnarray*}
   m^{\infty}-\frac{1}{N}\|\nabla \hat{u}\|_2^2+o(1)
    &  =  & \frac{1}{N}\|\nabla w_n\|_2^2\nonumber\\
    &  =  & I^{\infty}(w_n)-\frac{1}{N}\mathcal{P}^{\infty}(w_n)\nonumber\\
    & \ge & I^{\infty}\left({(w_n)}_{t_n}\right)-\frac{t_n^N}{N}\mathcal{P}^{\infty}(w_n)\nonumber\\
    & \ge & m^{\infty}-\frac{t_n^N}{N}\mathcal{P}^{\infty}(w_n)\ge m^{\infty},
 \end{eqnarray*}
 which implies $\mathcal{P}^{\infty}(\hat{u})\le 0$ due to $\|\nabla \hat{u}\|_2>0$. Since $\hat{u}\ne 0$ and
 $\mathcal{P}^{\infty}(\hat{u})\le 0$, in view of Lemma \ref{lem 2.6} and Corollary \ref{cor 2.8}, there exists
 $\hat{t}>0$ such that $\hat{u}_{\hat{t}}\in \mathcal{M}^{\infty}$.  From \eqref{Ii}, \eqref{Pi}, \eqref{G35}, \eqref{D71} and
 the weak semicontinuity of norm, one has
 \begin{eqnarray*}
   m^{\infty}
   &  =  & \lim_{n\to\infty} \left[I^{\infty}(\hat{u}_n)-\frac{1}{N}\mathcal{P}^{\infty}(\hat{u}_n)\right]\nonumber\\
   &  =  & \frac{1}{N}\lim_{n\to\infty} \|\nabla\hat{u}_n\|_2^2\ge \frac{1}{N}\|\nabla\hat{u}\|_2^2\nonumber\\
   &  =  & I^{\infty}(\hat{u})-\frac{1}{N}\mathcal{P}^{\infty}(\hat{u})\ge I^{\infty}\left({\hat{u}}_{\hat{t}}\right)
             -\frac{\hat{t}^N}{N}\mathcal{P}^{\infty}(\hat{u})\nonumber\\
   & \ge & m^{\infty}-\frac{\hat{t}^N}{N}\mathcal{P}^{\infty}(\hat{u})\ge m^{\infty},
 \end{eqnarray*}
 which implies
 \begin{equation*}
   \mathcal{P}^{\infty}(\hat{u})=0, \ \ \ \ I^{\infty}(\hat{u})=m^{\infty}.
 \end{equation*}
 \end{proof}

 \begin{lemma}\label{lem 2.13}
 Assume that {\rm(V1)-(V3)} and {\rm(F1)-(F3)} hold. If $\bar{u}\in \mathcal{M}$
 and $I(\bar{u})=m$, then $\bar{u}$ is a critical point of $I$.
\end{lemma}

\begin{proof} Assume that $I'(\bar{u})\ne 0$. Then there exist $\delta>0$ and $\varrho>0$ such that
 \begin{equation}\label{B81}
   \|u-\bar{u}\|\le 3\delta\Rightarrow \|I'(u)\|\ge \varrho.
 \end{equation}
 First, we prove that
 \begin{equation}\label{B82}
   \lim_{t\to 1}\left\|\bar{u}_t-\bar{u}\right\|=0.
 \end{equation}
 Arguing by contradiction, suppose that there exist $\varepsilon_0>0$ and a sequence $\{t_n\}$
 such that
 \begin{equation}\label{B83}
   \lim_{n\to \infty}t_n=1, \ \ \ \ \left\|\bar{u}_{t_n}-\bar{u}\right\|^2\ge \varepsilon_0.
 \end{equation}
 Since $\bar{u}\in H^1(\R^N)$, there exist $U\in \mathcal{C}_0(\R^N, \R^N)$ and $v\in \mathcal{C}_0(\R^N, \R)$ such that
 \begin{equation}\label{B84}
   \int_{\R^N}|\nabla\bar{u}-U|^2< \frac{\varepsilon_0}{20}, \ \ \ \ \int_{\R^N}|\bar{u}-v|^2< \frac{\varepsilon_0}{20}.
 \end{equation}
 From \eqref{B83} and \eqref{B84}, one has
 \begin{eqnarray}\label{B85}
   &     & \left\|\nabla\bar{u}_{t_n}-\nabla\bar{u}\right\|_2^2\nonumber\\
   &  =  & \int_{\R^N}\left|\nabla\bar{u}_{t_n}-\nabla\bar{u}\right|^2\mathrm{d}x \nonumber\\
   & \le & 2\int_{\R^N}\left|\nabla\bar{u}_{t_n}-U\right|^2\mathrm{d}x
            +2\int_{\R^N}|\nabla\bar{u}-U|^2\mathrm{d}x\nonumber\\
   &  =  & 2\int_{\R^N}\left|t_n^{-1}\nabla\bar{u}\left(t_n^{-1}x\right)-U(x)\right|^2\mathrm{d}x
            +2\int_{\R^N}|\nabla\bar{u}-U|^2\mathrm{d}x\nonumber\\
   & \le & 6t_n^{-2}\int_{\R^N}\left|U\left(t_n^{-1}x\right)-U(x)\right|^2\mathrm{d}x
            +6|t_n^{-1}-1|^2\int_{\R^N}|U|^2\mathrm{d}x+\frac{(1+3t_n^{N-2})\varepsilon_0}{10}\nonumber\\
   &  =  & \frac{2}{5}\varepsilon_0+o(1)
 \end{eqnarray}
 and
 \begin{eqnarray}\label{B86}
   \left\|\bar{u}_{t_n}-\bar{u}\right\|_2^2
   &  =  & \int_{\R^N}\left|\bar{u}_{t_n}-\bar{u}\right|^2\mathrm{d}x \nonumber\\
   & \le & 2\int_{\R^N}\left|\bar{u}_{t_n}-v\right|^2\mathrm{d}x +2\int_{\R^N}|\bar{u}-v|^2\mathrm{d}x\nonumber\\
   & \le & 4\int_{\R^N}\left|v\left(t_n^{-1}x\right)-v(x)\right|^2\mathrm{d}x+\frac{(1+2t_n^N)\varepsilon_0}{10}\nonumber\\
   &  =  & \frac{3}{10}\varepsilon_0+o(1).
 \end{eqnarray}
 Combining \eqref{B85} with \eqref{B86}, one has
 \begin{equation}\label{B87}
   \left\|\bar{u}_{t_n}-\bar{u}\right\|^2=\left\|\nabla(\bar{u}_{t_n})-\nabla\bar{u}\right\|_2^2
     +\left\|\bar{u}_{t_n}-\bar{u}\right\|_2^2\le \frac{7}{10}\varepsilon_0+o(1).
 \end{equation}
 \eqref{B87} contradicts with \eqref{B83}. Therefore, \eqref{B82} holds. Thus, there exists $\delta_1\in (0, 1/4)$ such that
 \begin{equation}\label{B89}
   |t-1|<\delta_1\Rightarrow \left\|\bar{u}_t-\bar{u}\right\|< \delta.
 \end{equation}
 In view of Lemma \ref{lem 2.2}, one has
 \begin{eqnarray}\label{B90}
   I\left(\bar{u}_t\right)
   & \le & I(\bar{u})-\frac{(1-\theta)\left[2-Nt^{N-2}+(N-2)t^N\right]}{2N}\|\nabla\bar{u}\|_2^2 \nonumber\\
   &  =  & m-\frac{(1-\theta)\mathfrak{g}(t)}{2N}\|\nabla\bar{u}\|_2^2, \ \ \ \ \forall \ t> 0.
 \end{eqnarray}
 It follows from \eqref{Jv}, \eqref{G38} and \eqref{G37} that there exist $T_1\in (0,1)$ and $T_2\in (1, \infty)$ such that
 \begin{equation}\label{T12}
   \mathcal{P}\left(\bar{u}_{T_1}\right)>0, \ \ \ \ \mathcal{P}\left(\bar{u}_{T_2}\right)<0.
 \end{equation}
 Let $\varepsilon:=\min\{(1-\theta)\mathfrak{g}(T_1)\|\nabla\bar{u}\|_2^2/5N, (1-\theta)\mathfrak{g}(T_2)\|\nabla\bar{u}\|_2^2/5N,
  1, \varrho\delta/8\}$ and $S:=B(\bar{u}, \delta)$. Then \cite[Lemma 2.3]{WM}
 yields a deformation $\eta\in \mathcal{C}([0, 1]\times H^1(\R^N), H^1(\R^N))$ such that

 \begin{enumerate}[i)]
  \item  $\eta(1, u)=u$ if $I(u)<m-2\varepsilon$ or $I(u)>m+2\varepsilon$;
  \item $\eta\left(1, I^{m+\varepsilon}\cap B(\bar{u}, \delta)\right)\subset I^{m-\varepsilon}$;
  \item $I(\eta(1, u))\le I(u), \ \forall \ u\in H^1(\R^N)$;
  \item $\eta(1, u)$ is a homeomorphism of $H^1(\R^N)$.
 \end{enumerate}
 By Corollary \ref{cor 2.4},  $I\left(\bar{u}_t\right)\le I(\bar{u})=m$ for $t> 0$, then it follows from \eqref{B89} and ii) that
 \begin{equation}\label{B91}
   I\left(\eta\left(1, \bar{u}_t\right)\right)\le m-\varepsilon, \ \ \ \ \forall \ t> 0, \ \ |t-1|< \delta_1.
 \end{equation}
 On the other hand, by iii) and \eqref{B90}, one has
 \begin{eqnarray}\label{B92}
   I\left(\eta\left(1, \bar{u}_t\right)\right)
   & \le & I\left(\bar{u}_t\right)\nonumber\\
   & \le & m-\frac{(1-\theta)\mathfrak{g}(t)}{2N}\|\nabla\bar{u}\|_2^2 \nonumber\\
   & \le & m-\frac{(1-\theta)\delta_2}{2N}\|\nabla\bar{u}\|_2^2, \ \ \ \ \forall \  t> 0, \ \ |t-1|\ge \delta_1,
 \end{eqnarray}
 where
 $$
   \delta_2:=\min\{\mathfrak{g}(1-\delta_1), \mathfrak{g}(1+\delta_1)\}>0.
 $$
 Combining \eqref{B91} with \eqref{B92}, we have
 \begin{equation}\label{B93}
   \max_{t\in [T_1, T_2]}I\left(\eta\left(1, \bar{u}_t\right)\right)<m.
 \end{equation}
 Define $\Psi_0(t):=\mathcal{P}\left(\eta\left(1, \bar{u}_t\right)\right)$ for $t> 0$.
 It follows from \eqref{B90} and i) that $\eta(1, \bar{u}_t)=\bar{u}_t$ for $t=T_1$ and $t=T_2$, which, together with \eqref{T12}, 
 implies
 $$
   \Psi_0(T_1)=\mathcal{P}\left(\bar{u}_{T_1}\right)>0, \ \ \ \ \Psi_0(T_2)=\mathcal{P}\left(\bar{u}_{T_1}\right)<0.
 $$
 Since $\Psi_0(t)$ is continuous on $(0, \infty)$, then  we have that $\eta\left(1, \bar{u}_t\right) \cap \mathcal{M}\ne \emptyset$ 
 for some $t_0\in [T_1, T_2]$, contradicting to the definition of $m$.
 \end{proof}

 \begin{proof}[Proof of Theorem  \ref{thm1.3}] In view of Lemmas \ref{lem 2.9}, \ref{lem 2.12} and \ref{lem 2.13}, there exists
 $\hat{u}\in \mathcal{M}^{\infty}$ such that
 $$
   I^{\infty}(\hat{u})=m^{\infty}=\inf_{u\in \Lambda}\max_{t > 0}I^{\infty}(u_t), \ \ \ \ (I^{\infty})'(\hat{u})=0.
 $$
 This shows that $\hat{u}$ is a ground state solution of Poho\u zaev type for \eqref{SE1}.
 \end{proof}

 \vskip6mm
 {\section{Ground state solutions for \eqref{SE}}}
 \setcounter{equation}{0}

 \vskip2mm
 \par
   In this section, we give the proof of Theorem \ref{thm1.2}. In the rest of two sections, we always assume that $V(x)\not \equiv V_{\infty}$
 in (V2) (if $V(x)\equiv V_{\infty}$, we recall that Theorem \ref{thm1.2} is contained in Theorem \ref{thm1.3}).

 \begin{lemma}\label{lem 3.1}
 Assume that {\rm(V1)-(V3)} and {\rm(F1)-(F3)} hold. Then $m^{\infty}\ge m$.
\end{lemma}

 \begin{proof} In view of Theorem \ref{thm1.3}, $I^{\infty}$ has a minimizer $u^{\infty}\ne 0$ on
 $\mathcal{M}^{\infty}$, i.e.
 \begin{equation}\label{P11}
   u^{\infty}\in \mathcal{M}^{\infty} \ \ \ \  \mbox{and} \ \ \ \ m^{\infty}=I^{\infty}(u^{\infty}).
 \end{equation}
 In view of Lemmas \ref{lem 2.6} and \ref{lem 2.7}, there exists $t_0>0$ such that $(u^{\infty})_{t_0}\in \mathcal{M}$. Thus,
 it follows from (V2), \eqref{IU}, \eqref{Ii}, \eqref{G35} and \eqref{P11} that
 $$
   m^{\infty}=I^{\infty}(u^{\infty})\ge I^{\infty}\left((u^{\infty})_{t_0}\right)
     \ge I\left((u^{\infty})_{t_0}\right)\ge m.
 $$
 Then $m^{\infty}\ge m$.
 \end{proof}

\begin{lemma}\label{lem 3.2}
 Assume that {\rm(V1)-(V3)} and {\rm(F1)-(F3)} hold. Then $m$ is achieved.
\end{lemma}

 \begin{proof} In view of Lemmas \ref{lem 2.6}, \ref{lem 2.7} and \ref{lem 2.11}, we have $\mathcal{M}\ne \emptyset$
 and $m>0$. Let $\{u_n\}\subset \mathcal{M}$ be such that $I(u_n)\rightarrow m$.
 Since $\mathcal{P}(u_n)=0$, then it follows from \eqref{G31} with $t\rightarrow 0$ that
 \begin{equation}\label{P24}
   m+o(1)= I(u_n)\ge \frac{1-\theta}{N}\|\nabla u_n\|_2^2.
 \end{equation}
 This shows that $\{\|\nabla u_n\|_2\}$ is bounded. Next, we prove that $\{\|u_n\|\}$ is also bounded.
 By (F1), (F2), \eqref{Jv}, \eqref{G36} and the Sobolev embedding theorem, one has
 \begin{eqnarray}\label{P25}
   \gamma_1\|u_n\|^2
   & \le & (N-2)\|\nabla u_n\|_2^2+\int_{{\R}^N}[NV(x)+\nabla V(x)\cdot x]u_n^2\mathrm{d}x\nonumber\\
   &  =  & 2N\int_{{\R}^N}F(u_n)\mathrm{d}x \nonumber\\
   & \le & \frac{\gamma_1}{2}\|u_n\|^2+C_4\|u_n\|_{2^*}^{2^*}\nonumber\\
   & \le & \frac{\gamma_1}{2}\|u_n\|^2+C_4S^{-2^*/2}\|\nabla u_n\|_2^{2^*}.
 \end{eqnarray}
 This shows that $\{u_n\}$ is bounded in $H^1(\R^N)$. Passing to a subsequence, we have
 $u_n\rightharpoonup \bar{u}$ in $H^1(\R^N)$. Then $u_n\rightarrow \bar{u}$ in $L_{\mathrm{loc}}^s(\R^N)$
 for $2\le s<2^*$ and $u_n\rightarrow \bar{u}$ a.e. in $\R^N$. There are two possible cases: i). $\bar{u}=0$
 and ii). $\bar{u}\ne 0$.

 \vskip2mm
 \par
   Case i). $\bar{u}=0$, i.e. $u_n\rightharpoonup 0$ in $H^1(\R^N)$. Then $u_n\rightarrow 0$ in $L_{\mathrm{loc}}^s(\R^N)$
 for $2\le s<2^*$ and $u_n\rightarrow 0$ a.e. in $\R^N$. By (V2) and \eqref{G37}, it is easy to show that
 \begin{equation}\label{P32}
   \lim_{n\to\infty}\int_{\R^N}[V_{\infty}-V(x)]u_n^2\mathrm{d}x=
   \lim_{n\to\infty}\int_{\R^N}\nabla V(x)\cdot xu_n^2\mathrm{d}x=0.
 \end{equation}
 From \eqref{IU}, \eqref{Ii}, \eqref{Pi}, \eqref{Jv} and \eqref{P32}, one can get
 \begin{equation}\label{P37}
   I^{\infty}(u_n)\rightarrow m, \ \ \ \ \mathcal{P}^{\infty}(u_n)\rightarrow 0.
 \end{equation}
 From Lemma \ref{lem 2.11} (i), \eqref{Pi} and \eqref{P37}, one has
 \begin{eqnarray}\label{P64}
   \min\{N-2, NV_{\infty}\}\rho_0^2
    & \le & \min\{N-2, NV_{\infty}\}\|u_n\|^2\nonumber\\
    & \le & (N-2)\|\nabla u_n\|_2^2+NV_{\infty}\|u_n\|_2^2\nonumber\\
    &  =  & 2N\int_{{\R}^N}F(u_n)\mathrm{d}x+o(1).
 \end{eqnarray}
 Using (F1), (F2), \eqref{P64} and Lions' concentration compactness principle \cite[Lemma 1.21]{WM}, we can prove that
 there exist $\delta>0$ and a sequence $\{y_n\}\subset \R^N$ such that $\int_{B_1(y_n)}|u_n|^2\mathrm{d}x> \delta$. Let
 $\hat{u}_n(x)=u_n(x+y_n)$. Then we have $\|\hat{u}_n\|=\|u_n\|$ and
 \begin{equation}\label{P65}
   \mathcal{P}^{\infty}(\hat{u}_n)= o(1), \ \ \ \ I^{\infty}(\hat{u}_n)\rightarrow m, \ \ \ \ \int_{B_1(0)}|\hat{u}_n|^2\mathrm{d}x> \delta.
 \end{equation}
 Therefore, there exists $\hat{u}\in H^1(\R^N)\setminus \{0\}$ such that, passing to a subsequence,
 \begin{equation}\label{P72}
 \left\{
   \begin{array}{ll}
     \hat{u}_n\rightharpoonup \hat{u}, & \mbox{in} \ H^1(\R^N); \\
     \hat{u}_n\rightarrow \hat{u}, & \mbox{in} \ L_{\mathrm{loc}}^s(\R^N), \ \forall \ s\in [1, 2^*);\\
     \hat{u}_n\rightarrow \hat{u}, & \mbox{a.e. on} \ \R^N.
   \end{array}
 \right.
 \end{equation}
 Let $w_n=\hat{u}_n-\hat{u}$. Then \eqref{P72} and Lemma \ref{lem 2.10} yield that \eqref{D73} and \eqref{D74} hold. Moreover,
 \begin{equation}\label{P75}
    \frac{1}{N}\|\nabla w_n\|_2^2=m-\frac{1}{N}\|\nabla \hat{u}\|_2^2+o(1),
      \ \ \ \ \mathcal{P}^{\infty}(w_n) = -\mathcal{P}^{\infty}(\hat{u})+o(1).
 \end{equation}
 If there exists a subsequence $\{w_{n_i}\}$ of $\{w_n\}$ such that $w_{n_i}=0$, then going to this subsequence, we have
 \begin{equation}\label{P76}
    I^{\infty}(\hat{u})=m, \ \ \ \ \mathcal{P}^{\infty}(\hat{u})=0.
 \end{equation}
 Next, we assume that $w_n\ne 0$. We claim that $\mathcal{P}^{\infty}(\hat{u})\le 0$.
 Otherwise, if $\mathcal{P}^{\infty}(\hat{u})>0$, then \eqref{P75} implies $\mathcal{P}^{\infty}(w_n) < 0$ for large $n$.
 In view of Lemma \ref{lem 2.6} and Corollary \ref{cor 2.8}, there exists $t_n>0$ such that $(w_n)_{t_n}\in \mathcal{M}^{\infty}$. From \eqref{Ii}, \eqref{Pi}, \eqref{G35} and \eqref{P75}, we obtain
 \begin{eqnarray*}
   m-\frac{1}{N}\|\nabla \hat{u}\|_2^2+o(1)
    &  =  & \frac{1}{N}\|\nabla w_n\|_2^2\nonumber\\
    &  =  & I^{\infty}(w_n)-\frac{1}{N}\mathcal{P}^{\infty}(w_n)\nonumber\\
    & \ge & I^{\infty}\left({(w_n)}_{t_n}\right)-\frac{t_n^N}{N}\mathcal{P}^{\infty}(w_n)\nonumber\\
    & \ge & m^{\infty}-\frac{t_n^N}{N}\mathcal{P}^{\infty}(w_n)\ge m^{\infty},
 \end{eqnarray*}
 which implies $\mathcal{P}^{\infty}(\hat{u})\le 0$ due to $\|\nabla \hat{u}\|_2>0$. Since $\hat{u}\ne 0$ and
 $\mathcal{P}^{\infty}(\hat{u})\le 0$, in view of Lemma \ref{lem 2.6} and Corollary \ref{cor 2.8}, there exists
 $\hat{t}>0$ such that $\hat{u}_{\hat{t}}\in \mathcal{M}^{\infty}$.  From \eqref{Ii}, \eqref{Pi}, \eqref{G35}, \eqref{P65} and
 the weak semicontinuity of norm, one has
 \begin{eqnarray*}
   m
   &  =  & \lim_{n\to\infty} \left[I^{\infty}(\hat{u}_n)-\frac{1}{N}\mathcal{P}^{\infty}(\hat{u}_n)\right]\nonumber\\
   &  =  & \frac{1}{N}\lim_{n\to\infty} \|\nabla\hat{u}_n\|_2^2\ge \frac{1}{N}\|\nabla\hat{u}\|_2^2\nonumber\\
   &  =  & I^{\infty}(\hat{u})-\frac{1}{N}\mathcal{P}^{\infty}(\hat{u})\ge I^{\infty}\left({\hat{u}}_{\hat{t}}\right)
             -\frac{\hat{t}^N}{N}\mathcal{P}^{\infty}(\hat{u})\nonumber\\
   & \ge & m^{\infty}-\frac{\hat{t}^N}{N}\mathcal{P}^{\infty}(\hat{u})\nonumber\\
   & \ge & m-\frac{\hat{t}^N}{N}\mathcal{P}^{\infty}(\hat{u})\ge m,
 \end{eqnarray*}
 which implies \eqref{P76} holds also. In view of Lemmas \ref{lem 2.6} and \ref{lem 2.7}, there exists $\tilde{t}>0$ such that
 $\hat{u}_{\tilde{t}}\in \mathcal{M}$, moreover, it follows from (V2), \eqref{IU}, \eqref{Ii}, \eqref{P76} and Corollary \ref{cor 2.3} that
 $$
   m\le I(\hat{u}_{\tilde{t}})\le I^{\infty}(\hat{u}_{\tilde{t}})\le I^{\infty}(\hat{u})=m.
 $$
 This shows that $m$ is achieved at $\hat{u}_{\tilde{t}}\in \mathcal{M}$.

 \vskip2mm
 \par
   Case ii). $\bar{u}\ne 0$. Let $v_n=u_n-\bar{u}$. Then Lemma \ref{lem 2.10} yields
 \begin{equation}\label{K71}
    I(u_n)=I(\bar{u})+I(v_n)+o(1)
 \end{equation}
 and
 \begin{equation}\label{K72}
   \mathcal{P}(u_n)=\mathcal{P}(\bar{u})+\mathcal{P}(v_n)+o(1).
 \end{equation}
 Set
 \begin{equation}\label{K73}
   \Psi(u)=\frac{1}{N}\|\nabla u\|_2^2-\frac{1}{2N}\int_{{\R}^N}(\nabla V(x), x)u^2\mathrm{d}x.
 \end{equation}
 Then it follows from  \eqref{V3} with $t=0$ and \eqref{G38} that
 \begin{equation}\label{K74}
   \Psi(u)\ge \frac{1-\theta}{N}\|\nabla u\|_2^2, \ \ \ \ \forall \ u\in H^1(\R^N).
 \end{equation}
 Since $I(u_n)\rightarrow m$ and $\mathcal{P}(u_n)=0$, then it follows from \eqref{IU}, \eqref{Jv}, \eqref{K71}, \eqref{K72} and \eqref{K73} that
 \begin{equation}\label{K75}
    \Psi(v_n)=m-\Psi(\bar{u})+o(1), \ \ \ \ \mathcal{P}(v_n) = -\mathcal{P}(\bar{u})+o(1).
 \end{equation}
 If there exists a subsequence $\{v_{n_i}\}$ of $\{v_n\}$ such that $v_{n_i}=0$, then going to this subsequence, we have
 \begin{equation}\label{K76}
    I(\bar{u})=m, \ \ \ \ \mathcal{P}(\bar{u})=0,
 \end{equation}
 which implies the conclusion of Lemma \ref{lem 3.2} holds. Next, we assume that $v_n\ne 0$. We claim that $\mathcal{P}(\bar{u})\le 0$.
 Otherwise $\mathcal{P}(\bar{u})>0$, then \eqref{K75} implies $\mathcal{P}(v_n) < 0$ for large $n$. In view of Lemmas \ref{lem 2.6} and
 \ref{lem 2.7}, there exists $t_n>0$ such that $(v_n)_{t_n}\in \mathcal{M}$.  From \eqref{IU}, \eqref{Jv}, \eqref{G31} and \eqref{K75},
  we obtain
 \begin{eqnarray*}
   m-\Psi(\bar{u})+o(1)
    &  =  & \Psi(v_n)\nonumber\\
    &  =  & I(v_n)-\frac{1}{N}\mathcal{P}(v_n)\nonumber\\
    & \ge & I\left({(v_n)}_{t_n}\right)-\frac{t_n^N}{N}\mathcal{P}(v_n)\nonumber\\
    & \ge & m-\frac{t_n^N}{N}\mathcal{P}(v_n)\ge m,
 \end{eqnarray*}
 which implies $\mathcal{P}(\bar{u})\le 0$ due to $\Psi(\bar{u})>0$. Since $\bar{u}\ne 0$ and $\mathcal{P}(\bar{u})\le 0$,
 in view of Lemmas \ref{lem 2.6} and \ref{lem 2.7}, there exists $\bar{t}>0$ such that $\bar{u}_{\bar{t}}\in \mathcal{M}$.
 From \eqref{IU}, \eqref{Jv}, \eqref{G31}, \eqref{K73}, \eqref{K74} and the weak semicontinuity of norm, one has
 \begin{eqnarray*}
   m
   &  =  & \lim_{n\to\infty} \left[I(u_n)-\frac{1}{N}\mathcal{P}(u_n)\right]\nonumber\\
   &  =  & \lim_{n\to\infty} \Psi(u_n)\ge \Psi(\bar{u})\nonumber\\
   &  =  & I(\bar{u})-\frac{1}{N}\mathcal{P}(\bar{u})\ge I\left({\bar{u}}_{\bar{t}}\right)
             -\frac{\bar{t}^N}{N}\mathcal{P}(\bar{u})\nonumber\\
   & \ge & m-\frac{\bar{t}^N}{N}\mathcal{P}(\bar{u})\ge m,
 \end{eqnarray*}
 which implies \eqref{K76} also holds.
 \end{proof}

 \begin{proof}[Proof of Theorem \ref{thm1.2}] In view of Lemmas \ref{lem 2.9}, \ref{lem 2.13} and \ref{lem 3.2},
 there exists $\bar{u}\in \mathcal{M}$ such that
 $$
   I(\bar{u})=m=\inf_{u\in \Lambda}\max_{t > 0}I(u_t), \ \ \ \ I'(\bar{u})=0.
 $$
 This shows that $\bar{u}$ is a ground state solution of Poho\u zaev type for \eqref{SE}.
 \end{proof}

 \vskip6mm
 {\section{The least energy solutions for \eqref{SE}}}
 \setcounter{equation}{0}

 \vskip2mm
 \par
   In this section, we give the proof of Theorem \ref{thm1.5}.

 \begin{proposition}\label{pro 4.1}{\rm\cite{JTo}}
 Let $X$ be a Banach space and let $J\subset \R^+$ be an interval, and
 $$
   \Phi_{\lambda}(u)=A(u)-\lambda B(u), \ \ \ \ \forall \ \lambda\in J,
 $$
 be a family of $\mathcal{C}^1$-functional on $X$ such that
 \begin{enumerate}[{\rm i)}]
  \item either $A(u)\to +\infty$ or $B(u)\to +\infty$, as $\|u\| \to \infty$;
  \item $B$ maps every bounded set of $X$ into a set of $\R$ bounded below;
  \item there are two points $v_1, v_2$ in $X$ such that
 \begin{equation}\label{cm}
   \tilde{c}_{\lambda}:=\inf_{\gamma\in \tilde{\Gamma}}\max_{t\in [0, 1]}\Phi_{\lambda}(\gamma(t))>\max\{\Phi_{\lambda}(v_1), \Phi_{\lambda}(v_2)\},
 \end{equation}
 \end{enumerate}
 where
 $$
   \tilde{\Gamma}=\left\{\gamma\in \mathcal{C}([0, 1], X): \gamma(0)=v_1, \gamma(1)=v_2\right\}.
 $$
 Then, for almost every $\lambda\in J$, there exists a sequence $\{u_n(\lambda)\}$ such that
 \begin{enumerate}[{\rm i)}]
  \item $\{u_n(\lambda)\}$ is bounded in $X$;
  \item $\Phi_{\lambda}(u_n(\lambda))\rightarrow c_{\lambda}$;
  \item $\Phi_{\lambda}'(u_n(\lambda))\rightarrow 0$ in $X^*$, where $X^*$ is the dual of $X$.
 \end{enumerate}
\end{proposition}

\begin{lemma}\label{lem 4.2}{\rm\cite{JT2}}
 Assume that {\rm(V1), (V2), (F1) and (F2)}  hold.  Let $u$ be a critical point of
 $I_{\lambda}$ in $H^1(\R^N)$, then we have the following Poho\u zaev type identity
 \begin{eqnarray}\label{Pl}
  \mathcal{P}_{\lambda}(u)
    & :=  & \frac{N-2}{2}\|\nabla u\|_2^2+\frac{1}{2}\int_{\R^N}\left[NV(x)+\nabla V(x)\cdot x\right]u^2\mathrm{d}x\nonumber\\
    &     & \ \  -N\lambda\int_{\R^N}F(u)\mathrm{d}x=0.
 \end{eqnarray}
\end{lemma}

 \vskip4mm
 \par
    By Corollary \ref{cor 2.3}, we have the following lemma.

\begin{lemma}\label{lem 4.3}
 Assume that {\rm(F1)} and {\rm(F2)}  hold. Then
 \begin{eqnarray}\label{G45}
   I_{\lambda}^{\infty}(u)
    & \ge & I_{\lambda}^{\infty}\left(u_t\right)+\frac{1-t^N}{N}\mathcal{P}_{\lambda}^{\infty}(u)
              +\frac{2-Nt^{N-2}+(N-2)t^N}{2N}\|\nabla u\|_2^2,\nonumber\\
    &     &   \ \ \ \ \forall \ u\in H^1(\R^N), \ \ t> 0, \ \ \lambda\ge 0.
 \end{eqnarray}
\end{lemma}

 \vskip4mm
 \par
    In view of Theorem \ref{thm1.3}, $I_1^{\infty}=I^{\infty}$ has a minimizer $u_1^{\infty}\ne 0$ on $\mathcal{M}_1^{\infty}
 =\mathcal{M}^{\infty}$,  i.e.
 \begin{equation}\label{G47}
   u_1^{\infty}\in \mathcal{M}_1^{\infty}, \ \ \ \ (I_1^{\infty})'(u_1^{\infty})=0
     \ \ \ \  \mbox{and} \ \ \ \ m_1^{\infty}=I_1^{\infty}(u_1^{\infty}),
 \end{equation}
 where $m_{\lambda}^{\infty}$ is defined by \eqref{cm1}. Since \eqref{SE1} is autonomous, $V\in \mathcal{C}(\R^N, \R)$ and
 $V(x)\le V_{\infty}$ but $V(x)\not\equiv V_{\infty}$, then there exist $\bar{x}\in \R^N$ and $\bar{r}>0$ such that
 \begin{equation}\label{G48}
    V_{\infty}-V(x)>0, \ \ |u_1^{\infty}(x)|>0\ \ \ \ a.e. \ |x-\bar{x}|\le \bar{r}.
 \end{equation}

 \begin{lemma}\label{lem 4.4}
 Assume that {\rm(V1), (V2)} and {\rm(F1)-(F3)} hold. Then
 \begin{enumerate}[{\rm(i)}]
  \item there exists $T>0$ independent of $\lambda$ such that $I_{\lambda}\left((u_1^{\infty})_{T}\right)<0$
 for all $\lambda\in [0.5, 1]$;
  \item there exists a positive constant $\kappa_0 $ independent of $\lambda$ such that for all $\lambda\in [0.5, 1]$,
 \begin{equation*}
   c_{\lambda}:=\inf_{\gamma\in \Gamma}\max_{t\in [0, 1]}I_{\lambda}(\gamma(t))\ge \kappa_0
     >\max\left\{I_{\lambda}(0), I_{\lambda}\left((u_1^{\infty})_{T}\right)\right\},
 \end{equation*}
 where
 $$
   \Gamma=\left\{\gamma\in \mathcal{C}([0, 1], H^1(\R^N)): \gamma(0)=0, \gamma(1)=(u_1^{\infty})_{T}\right\};
 $$
 \item $c_{\lambda}$ is bounded for $\lambda\in [0.5, 1]$;
 \item $m_{\lambda}^{\infty}$ is non-increasing on $\lambda\in [0.5, 1]$;
 \item $\limsup_{\lambda\to \lambda_0}c_{\lambda}\le c_{\lambda_0}$ for $\lambda_0\in(0.5, 1]$.
 \end{enumerate}
\end{lemma}

 \vskip4mm
 \par
   Since $m_{\lambda}^{\infty}=I_{\lambda}^{\infty}(u_{\lambda}^{\infty})$ and $\int_{\R^N}F(u_{\lambda}^{\infty})\mathrm{d}x>0$,
 then the proofs of (i)-(iv) in Lemma \ref{lem 4.4} is standard, (v) can be proved similar to \cite[Lemma 2.3]{Je}, so we omit it.

 \begin{lemma}\label{lem 4.5}
 Assume that {\rm(V1), (V2)} and {\rm(F1)-(F3)} hold. Then
 there exists $\bar{\lambda}\in [1/2, 1)$ such that $c_{\lambda}<m_{\lambda}^{\infty}$ for $\lambda\in (\bar{\lambda}, 1]$.
\end{lemma}

\begin{proof} It is easy to see that $I_{\lambda}\left((u_1^{\infty})_t\right)$ is continuous on $t\in (0, \infty)$.
 Hence for any $\lambda\in [1/2, 1]$,  we can choose $t_{\lambda}\in (0, T)$ such that $I_{\lambda}
 \left((u_1^{\infty})_{t_{\lambda}}\right) =\max_{t\in [0,T]}I_{\lambda}\left((u_1^{\infty})_t\right)$. Setting
 $$
   \gamma_0(t)=\left\{\begin{array}{ll}
    (u_1^{\infty})_{(tT)}, \ \ &\mbox{for} \ t>0,\\
    0, \ \ & \mbox{for} \ t=0.
    \end{array}\right.
 $$
 Then $\gamma_0\in \Gamma$ defined by Lemma \ref{lem 4.4} (ii). Moreover
 \begin{equation}\label{G50}
   I_{\lambda} \left((u_1^{\infty})_{t_{\lambda}}\right)=\max_{t\in [0,1]}I_{\lambda}\left(\gamma_0(t)\right)
       \ge c_{\lambda}.
 \end{equation}
 Let
 \begin{equation}\label{G56}
   \zeta_0:=\min\{3\bar{r}/8(1+|\bar{x}|), 1/4\}.
 \end{equation}
 Then it follows from \eqref{G48} and \eqref{G56} that
 \begin{equation}\label{G58}
   |x-\bar{x}|\le \frac{\bar{r}}{2} \ \ \mbox{and} \ \ s\in [1-\zeta_0, 1+\zeta_0] \Rightarrow |sx-\bar{x}|\le \bar{r}.
 \end{equation}
 Since $\mathcal{P}^{\infty}(u_1^{\infty})=0$, then $\int_{\R^N}F(u_1^{\infty})\mathrm{d}x>0$. Let
 \begin{eqnarray}\label{G59}
   \bar{\lambda}
    & :=  & \max\left\{\frac{1}{2}, 1-\frac{(1-\zeta_0)^N\min_{s\in [1-\zeta_0, 1+\zeta_0]}\int_{\R^N}\left[V_{\infty}-
              V(sx)\right]|u_1^{\infty}|^2\mathrm{d}x}{T^{N}\int_{\R^N}F(u_1^{\infty})\mathrm{d}x},\right.\nonumber\\
    &     & \left.1-\frac{\min\{\mathfrak{g}(1-\zeta_0),\mathfrak{g}(1+\zeta_0)\}\|\nabla u_1^{\infty}\|_2^2}{N
             T^{N}\int_{\R^N}F(u_1^{\infty})\mathrm{d}x}\right\}. \ \ \ \
 \end{eqnarray}
 Then it follows from \eqref{B20}, \eqref{G48} and \eqref{G58} that $1/2\le \bar{\lambda}<1$. We have two cases to distinguish:

 \vskip2mm
 \par
   Case i). $t_{\lambda}\in [1-\zeta_0, 1+\zeta_0]$. From \eqref{Ilu}, \eqref{Il}, \eqref{G45}-\eqref{G50}, \eqref{G58}, \eqref{G59} and
 Lemma \ref{lem 4.4} (iv), we have
 \begin{eqnarray*}
   m_{\lambda}^{\infty}
    & \ge & m_1^{\infty}=I_1^{\infty}(u_1^{\infty})\ge I_1^{\infty}\left((u_1^{\infty})_{t_{\lambda}}\right)\nonumber\\
    &  =  & I_{\lambda}\left((u_1^{\infty})_{t_{\lambda}}\right)
              -\frac{(1-\lambda)t_{\lambda}^{N}}{2}\int_{\R^N}F(u_1^{\infty})\mathrm{d}x  +\frac{t_{\lambda}^N}{2}\int_{\R^N}[V_{\infty}-V(t_{\lambda}x)]|u_1^{\infty}|^2\mathrm{d}x\nonumber\\
    & \ge & c_{\lambda} -\frac{(1-\lambda)T_0^{N}}{2}\int_{\R^N}F(u_1^{\infty})\mathrm{d}x\nonumber\\
    &     & \ \     +\frac{(1-\zeta_0)^N}{2}\min_{s\in [1-\zeta_0, 1+\zeta_0]}
                \int_{\R^N}\left[V_{\infty}-V(sx)\right]|u_1^{\infty}|^2\mathrm{d}x\nonumber\\
    &  >  & c_{\lambda}, \ \ \ \ \forall \ \lambda\in (\bar{\lambda}, 1].
 \end{eqnarray*}
 \par
   Case ii). $t_{\lambda}\in (0, 1-\zeta_0)\cup (1+\zeta_0, T)$. From (V2), \eqref{Ilu}, \eqref{Il}, \eqref{B20}, \eqref{G45}, \eqref{G47},
 \eqref{G50}, \eqref{G59} and Lemma \ref{lem 4.4} (iv), we have
 \begin{eqnarray*}
   m_{\lambda}^{\infty}
    & \ge & m_1^{\infty}=I_1^{\infty}(u_1^{\infty})\ge I_1^{\infty}\left((u_1^{\infty})_{t_{\lambda}}\right)
             +\frac{\mathfrak{g}(t_{\lambda})\|\nabla u_1^{\infty}\|_2^2}{2N}\nonumber\\
    &  =  & I_{\lambda}\left((u_1^{\infty})_{t_{\lambda}}\right)
              -\frac{(1-\lambda)t_{\lambda}^{N}}{2}\int_{\R^N}F(u_1^{\infty})\mathrm{d}x\nonumber\\
    &     & \ \  +\frac{t_{\lambda}^N}{2}\int_{\R^N}[V_{\infty}-V(t_{\lambda}x)]|u_1^{\infty}|^2\mathrm{d}x
              +\frac{\mathfrak{g}(t_{\lambda})\|\nabla u_1^{\infty}\|_2^2}{2N}\nonumber\\
    & \ge & c_{\lambda} -\frac{(1-\lambda)T^{N}}{2}\int_{\R^N}F(u_1^{\infty})\mathrm{d}x\nonumber\\
    &     & \ \  +\frac{\min\{\mathfrak{g}(1-\zeta_0),\mathfrak{g}(1+\zeta_0)\}\|\nabla u_1^{\infty}\|_2^2}{2N}\nonumber\\
    &  >  & c_{\lambda}, \ \ \ \ \forall \ \lambda\in (\bar{\lambda}, 1].
 \end{eqnarray*}
 In both cases, we obtain that $c_{\lambda}<m_{\lambda}^{\infty}$ for $\lambda\in (\bar{\lambda}, 1]$.
 \end{proof}

 \begin{lemma}\label{lem 4.6}{\rm\cite{JT2}}
 Assume that {\rm(V1), (V2)} and {\rm(F1)-(F3)} hold. Let $\{u_n\}$ be a bounded (PS)
 sequence for $I_{\lambda}$, for $\lambda\in  [1/2, 1]$. Then there exists a subsequence of $\{u_n\}$, still denoted by
 $\{u_n\}$, an integer $l\in \N \cup \{0\}$, a sequence $\{y_n^k\}$ and $w^k\in H^1(\R^3)$ for $1\le k\le l$, such that
 \begin{enumerate}[{\rm(i)}]
  \item $u_n\rightharpoonup u_0$ with $I_{\lambda}'(u_0)=0$;
  \item $w^k\ne 0$ and $(I_{\lambda}^{\infty})'(w^k)=0$ for $1\le k\le l$;
  \item $\left\|u_n-u_0-\sum_{k=1}^lw^k(\cdot+y_n^k)\right\|\rightarrow 0$;
  \item $I_{\lambda}(u_n)\rightarrow I_{\lambda}(u_0)+\sum_{i=1}^{l}I_{\lambda}^{\infty}(w^i)$;
 \end{enumerate}
 where we agree that in the case $l = 0$ the above holds without $w^k$.
\end{lemma}

 \begin{lemma}\label{lem 4.7}
 Assume that {\rm(V1), (V2), (V4)} and {\rm(F1)-(F3)} hold.Then for almost every
 $\lambda\in (\bar{\lambda},1]$, there exists $u_{\lambda}\in H^1(\R^N)\setminus \{0\}$ such that
 \begin{equation}\label{D31}
   I_{\lambda}'(u_{\lambda})=0, \ \ \ \ I_{\lambda}(u_{\lambda}) = c_{\lambda}.
 \end{equation}
\end{lemma}

\begin{proof}  Under (V1), (V2) and (F1)-(F3), Lemma \ref{lem 4.4} implies that $I_{\lambda}(u)$ satisfies the assumptions of Proposition \ref{pro 4.1} with $X=H^1(\R^N)$, $J=[\bar{\lambda},1]$ and $\Phi_{\lambda}=I_{\lambda}$. So for almost every $\lambda\in (\bar{\lambda},1]$, there exists a bounded sequence $\{u_n(\lambda)\} \subset H^1(\R^N)$ (for simplicity, we denote the sequence by $\{u_n\}$ instead of $\{u_n(\lambda)\}$)
 such that
 \begin{equation}\label{PS}
   I_{\lambda}(u_n)\rightarrow c_{\lambda}>0, \ \ \ \ I_{\lambda}'(u_n) \rightarrow 0.
 \end{equation}
 By Lemmas \ref{lem 4.2} and \ref{lem 4.6}, there exist a subsequence of $\{u_n\}$, still denoted by $\{u_n\}$, and $u_{\lambda}\in H^1(\R^N)$,
 an integer $l\in \N\cup \{0\}$, and $w^1, \ldots, w^l\in H^1(\R^N)\setminus \{0\}$ such that
 \begin{equation}\label{un1}
   u_n\rightharpoonup u_{\lambda}\ \ \mbox {in} \  H^1(\R^N), \ \ \ \  I_{\lambda}'(u_{\lambda})=0,
 \end{equation}
 \begin{equation}\label{un2}
   (I_{\lambda}^{\infty})'(w^k)=0, \ \ \ \ I_{\lambda}^{\infty}(w^k)\ge m_{\lambda}^{\infty},\ \ \ \ 1\le k\le l
 \end{equation}
 and
 \begin{equation}\label{Ab1}
    c_{\lambda}= I_{\lambda}(u_{\lambda})+\sum_{k=1}^{l}I_{\lambda}^{\infty}(w^k).
 \end{equation}

 \par
 Since $I_{\lambda}'(u_{\lambda})=0$, then it follows from Lemma \ref{lem 4.2} that
 \begin{eqnarray}\label{Plt}
  \mathcal{P}_{\lambda}(u_{\lambda})
    &  =  & \frac{N-2}{2}\|\nabla u_{\lambda}\|_2^2+\frac{1}{2}\int_{\R^N}\left[NV(x)+\nabla V(x)\cdot x\right]u_{\lambda}^2\mathrm{d}x\nonumber\\
    &     & \ \  -N\lambda\int_{\R^N}F(u_{\lambda})\mathrm{d}x=0.
 \end{eqnarray}
 Since $\|u_n\|\nrightarrow 0$, we deduce from \eqref{un2} and \eqref{Ab1} that if $u_{\lambda}=0$ then $l\ge 1$ and
 \begin{eqnarray*}
   c_{\lambda} =  I_{\lambda}(u_{\lambda})+\sum_{k=1}^{l}I_{\lambda}^{\infty}(w^k)
    \ge  m_{\lambda}^{\infty},
 \end{eqnarray*}
 which contradicts with Lemma \ref{lem 4.5}. Thus $u_{\lambda}\ne 0$.
 It follows from \eqref{Ilu}, \eqref{G38}, \eqref{Plt} and (V4) that
 \begin{eqnarray}\label{D33}
   I_{\lambda}(u_{\lambda})
     &  =  & I_{\lambda}(u_{\lambda})-\frac{1}{N}\mathcal{P}_{\lambda}(u_{\lambda})\nonumber\\
     &  =  & \frac{1}{N}\|\nabla u_{\lambda}\|_2^2-\frac{1}{2N}\int_{{\R}^N}(\nabla V(x), x)u_{\lambda}^2\mathrm{d}x
               \ge \frac{1-\theta}{N}\|\nabla u_{\lambda}\|_2^2>0.
 \end{eqnarray}
 From \eqref{Ab1} and \eqref{D33}, one has
 \begin{eqnarray}\label{D34}
   c_{\lambda} =  I_{\lambda}(u_{\lambda})+\sum_{k=1}^{l}I_{\lambda}^{\infty}(w^k)
    >  lm_{\lambda}^{\infty}.
 \end{eqnarray}
 By Lemma \ref{lem 4.5}, we have $c_{\lambda}<m_{\lambda}^{\infty}$ for $\lambda\in (\bar{\lambda}, 1]$, which, together with
 \eqref{D34}, implies that $l=0$ and $I_{\lambda}(u_{\lambda}) = c_{\lambda}$.
 \end{proof}

 \begin{lemma}\label{lem 4.8}
 Assume that {\rm(V1), (V2), (V4)} and {\rm(F1)-(F3)} hold. Then there exists
 $\bar{u}\in H^1(\R^N)\setminus \{0\}$ such that
 \begin{equation}\label{Q05}
   I'(\bar{u})=0, \ \ \ \ 0< I(\bar{u}) < c_1.
 \end{equation}
\end{lemma}

\begin{proof} In view of Lemma \ref{lem 4.7}, there exist two sequences $\{\lambda_n\}\subset [\bar{\lambda}, 1]$
 and $\{u_{\lambda_n}\}\subset H^1(\R^N)\setminus \{0\}$, denoted by $\{u_n\}$, such that
 \begin{equation}\label{Q00}
   \lambda_n\rightarrow 1, \ \ \ \ c_{\lambda_n}\rightarrow c_*, \ \ \ \ I_{\lambda_n}'(u_n)=0, \ \ \ \ 0<I_{\lambda_n}(u_n) \le c_{\lambda_n}.
 \end{equation}
 Then it follows from Lemma \ref{lem 4.2} that
 \begin{eqnarray}\label{Pln}
  \mathcal{P}_{\lambda_n}(u_n)
    & :=  & \frac{N-2}{2}\|\nabla u_n\|_2^2+\frac{1}{2}\int_{\R^N}\left[NV(x)+\nabla V(x)\cdot x\right]u_n^2\mathrm{d}x\nonumber\\
    &     & \ \  -N\lambda_n\int_{\R^N}F(u_n)\mathrm{d}x=0.
 \end{eqnarray}
 From  (V4), \eqref{Ilu}, \eqref{G38}, \eqref{Q00}, \eqref{Pln} and Lemma \ref{lem 4.4} (iii), one has
 \begin{eqnarray}\label{Q01}
   C_5
    & \ge & c_{\lambda_n}=I_{\lambda_n}(u_n)-\frac{1}{N}\mathcal{P}_{\lambda_n}(u_n)\nonumber\\
    &  =  & \frac{1}{N}\|\nabla u_n\|_2^2-\frac{1}{2N}\int_{\R^N}\nabla V(x)\cdot xu_n^2\mathrm{d}x\nonumber\\
    & \ge & \frac{1-\theta}{N}\|\nabla u_n\|_2^2.
 \end{eqnarray}
 This shows that $\{\|\nabla u_n\|_2\}$ is bounded.  Next, we demonstrate that $\{u_n\}$ is bounded in $H^1(\R^N)$.
 According to (V1) and (V2), it is easy to show that there exists a constant $\gamma_3>0$ such that
 \begin{equation}\label{Q02}
    \int_{\R^N}\left[|\nabla u|^2+V(x)u^2\right]\mathrm{d}x\ge \gamma_3\|u\|^2, \ \ \ \ \forall \ u\in H^1(\R^N).
 \end{equation}
 From (F1), (F2), \eqref{Ilu}, \eqref{Q00}, \eqref{Q01}, \eqref{Q02}, Lemma \ref{lem 4.4} (iii) and the Sobolev embedding theorem, we have
 \begin{eqnarray*}\label{Q03}
   \gamma_3\|u_n\|^2
    & \le & \int_{\R^N}\left[|\nabla u_n|^2+V(x)u_n^2\right]\mathrm{d}x\\
    &  =  & 2c_{\lambda_n}+2\lambda_n\int_{\R^N}F(u_n)\mathrm{d}x\\
    & \le & 2C_5+\frac{\gamma_3}{2}\|u_n\|^2+C_6\|u_n\|_{2^*}^{2^*}\nonumber\\
   & \le & 2C_5+\frac{\gamma_3}{2}\|u_n\|^2+C_6S^{-2^*/2}\|\nabla u_n\|_2^{2^*}.
 \end{eqnarray*}
 Hence, $\{u_n\}$ is bounded in $H^1(\R^N)$. In view of Lemma \ref{lem 4.4} (v), we have $\lim_{n\to\infty}c_{\lambda_n}=c_*\le c_1$. Hence,
 it follows from \eqref{IU}, \eqref{Ilu} and \eqref{Q00} that
 \begin{equation}\label{Q04}
   I(u_n)\rightarrow  c_*, \ \ \ \  I'(u_n)\rightarrow 0.
 \end{equation}
 This shows that $\{u_n\}$ satisfy \eqref{PS} with $c_{\lambda}=c_*$. In view of the proof of Lemma \ref{lem 4.7}, we can show that there
 exists $\bar{u}\in H^1(\R^N)\setminus \{0\}$ such that \eqref{Q05} holds.
 \end{proof}

 \begin{proof}[Proof of Theorem  \ref{thm1.5}]  Let
 $$
   \mathcal{K}:=\left\{u\in H^1(\R^N)\setminus \{0\} : I'(u)=0\right\}, \ \ \ \ \hat{m}:=\inf_{u\in\mathcal{K}}I(u).
 $$
 Then Lemma \ref{lem 4.8} shows that $\mathcal{K}\ne \emptyset$ and $\hat{m}\le c_1$. For any $u\in \mathcal{K}$, Lemma \ref{lem 4.2}
 implies $\mathcal{P}(u)=\mathcal{P}_1(u)=0$. Hence it follows from \eqref{D33} that $I(u)=I_1(u)>0$, and so $\hat{m}\ge 0$.
 Let $\{u_n\}\subset \mathcal{K}$ such that
 \begin{equation}\label{Z01}
   I'(u_n)=0, \ \ \ \ I(u_n) \rightarrow \hat{m}.
 \end{equation}
 In view of Lemma \ref{lem 4.5}, $\hat{m}\le c_1<m_1^{\infty}$. By a similar argument as in the proof of Lemma \ref{lem 4.7}, we can prove
 that there exists $\bar{u}\in H^1(\R^N)\setminus \{0\}$ such that
 \begin{equation}\label{Z02}
   I'(\bar{u})=0, \ \ \ \ I(\bar{u}) = \hat{m}.
 \end{equation}
 This shows that $\bar{u}$ is a nontrivial least energy solution of \eqref{SE}.
 \end{proof}

\bibliographystyle{plain}

\bibliography{TXH1706}

\begin{thebibliography}{10}

\bibitem{BL}
H.~Berestycki and P.~L. Lions.
\newblock Nonlinear scalar field equations, {I}. {E}xistence of a ground state.
\newblock {\em Arch. Rational Mech. Anal.}, 82:313--345, 1983.

\bibitem{Je}
L.~Jeanjean.
\newblock On the existence of bounded {P}alais-{S}male sequences and
  application to a {L}andesman- {L}azer-type problem set on {$\mathbb R^N$}.
\newblock {\em Proc. Roy. Soc. Edinburgh Sect. A}, 129:787--809, 1999.

\bibitem{JT}
L.~Jeanjean and K.~Tanka.
\newblock A remark on least energy solutions in {$\mathbb R^N$}.
\newblock {\em Proc. Amer. Math. Soc.}, 131:2399--2408, 2003.

\bibitem{JT2}
L.~Jeanjean and K.~Tanka.
\newblock A positive solution for a nonlinear {S}chr\"odinger equation on
  {$\mathbb R^N$}.
\newblock {\em Indiana Univ. Math. J.}, 54:443--464, 2005.

\bibitem{JTo}
L.~Jeanjean and J.~F. Toland.
\newblock Bounded {P}alais-{S}male mountain-pass sequences.
\newblock {\em C. R. Acad. Sci. Paris S\'er. I Math.}, 327:23--28, 1998.

\bibitem{Po}
S.~I. Poho{\u z}aev.
\newblock Eigenfunctions of the equation $au + 2f(u) = 0$.
\newblock {\em Sov. Math. Doklady}, 5:1408--1411, 1965.

\bibitem{Ra}
P.~H. Rabinowitz.
\newblock On a class of nonlinear {S}chr\"odinger equations.
\newblock {\em Z. Angew. Math. Phys.}, 43:270--291, 1992.

\bibitem{Sh}
J.~Shatah.
\newblock Unstable ground state of nonlinear {K}lein-{G}ordon equations.
\newblock {\em Trans. Amer. Math. Soc.}, 290:701--710, 1985.

\bibitem{St}
M.~Struwe.
\newblock {\em Variational methods}.
\newblock Results in Mathematics and Related Areas, 3. Springer-Verlag, Berlin,
  1996.

\bibitem{Ta}
X.~H. Tang.
\newblock Non-{N}ehari manifold method for asymptotically periodic
  {S}chr\"odinger equations.
\newblock {\em Sci. China Math.}, 58:715--728, 2015.

\bibitem{TC2}
X.~H. Tang and S.~T. Chen.
\newblock Ground state solutions of {N}ehari-{P}oho\u zaev type for
  {K}irchhoff-type problems with general potentials.
\newblock {\em Calc. Var. Partial Differential Equations}, 56:110--134.

\bibitem{TC1}
X.~H. Tang and S.~T. Chen.
\newblock Ground state solutions of {N}ehari-{P}oho\u zaev type for
  {S}chr\"odinger-{P}oisson problems with general potentials.
\newblock {\em Disc. Contin. Dyn. Syst.}, 37:4973--5002, 2017.

\bibitem{WM}
M.~Willem.
\newblock {\em Minimax theorems}.
\newblock Progress in Nonlinear Differential Equations and their Applications,
  24. Birkh\"auser Boston Inc., Boston, MA, 1996.

\end{thebibliography}

\end{document}